\documentclass[a4paper,oneside,11pt]{article}
\usepackage{a4wide, dsfont}
\usepackage{natbib}
\usepackage{tikz}
\usepackage{epsfig}
\usepackage{graphicx}
\usepackage[latin1]{inputenc}
\usepackage{amsmath}%
\usepackage{subfigure}
\usepackage{url}
\usepackage{amssymb}%
\usepackage{amsthm}%
\theoremstyle{plain}
\newtheorem{theorem}{Theorem}%
\newtheorem{cor}[theorem]{Corollary}%
\theoremstyle{definition}
\newtheorem{remark}{Remark}%
\newtheorem*{remark*}{Remark}%
\usepackage{diagbox}
\usepackage{threeparttable}
\DeclareMathOperator{\diag}{diag}

\newcommand\f[1]{\boldsymbol{#1}}
\newcommand{\fpi}{\boldsymbol{\pi}}

\newcommand{\cond}{\,|\,}

\renewcommand{\theta}{\vartheta}
\newcommand{\cC}{\mathcal{C}}

\newcommand{\bx}{\ensuremath{\boldsymbol{x}}}%
\newcommand{\bp}{\ensuremath{\boldsymbol{p}}}%
\newcommand{\bz}{\ensuremath{\boldsymbol{z}}}%
\newcommand{\bX}{\ensuremath{\boldsymbol{X}}}%
\newcommand{\bZ}{\ensuremath{\boldsymbol{Z}}}%

\newcommand{\N}{\ensuremath{\mathbb{N}}}%
\newcommand{\cP}{\ensuremath{\mathcal{P}}}%
\newcommand{\cX}{\ensuremath{\mathcal{X}}}%
\renewcommand\theta{\vartheta}

\begin{document}
%
%
{\title{Nonparametric identification in the dynamic stochastic block model}
	\author{Ann-Kristin Becker  \\
		\small{Institute of Bioinformatics}  \\
		\small{University Medicine Greifswald} \\
		\small{ann-kristin.becker@uni-greifswald.de} 
			\and
	Hajo Holzmann\footnote{Corresponding author. Prof.~Dr.~Hajo Holzmann,  Faculty of Mathematics and Computer Science, Philipps-Universität Marburg, Hans-Meerweinstr., 35043 Marburg, Germany} \\
	\small{Fachbereich Mathematik und Informatik}  \\
	\small{Philipps-Universität Marburg} \\
	\small{holzmann@mathematik.uni-marburg.de}}
	\maketitle
}
%
%

%
\normalsize
\begin{abstract}
We show nonparametric identification of the parameters in the dynamic stochastic block model as recently introduced in \citet{Matias2016} in case of binary, finitely weighted and general edge states. We formulate conditions on the true parameters which guarantee actual point identification instead of mere generic identification, and which also lead to novel conclusions in the static case. In particular, our results justify in terms of identification the applications in \citet{Matias2016} to finitely weighted edges with three edge states. We also give numerical illustrations via the variational EM algorithm in simulation settings covered by our identification analysis. 
\end{abstract}
%

{\sl Key words: dynamic network, dynamic random graph, stochastic block model, community detection }  \\*[0.5cm]


%
%
%
\section{Introduction}
Nowadays network data occurs frequently, e.g.~in the form of social, communication, ecological or protein networks, and much effort has recently been devoted to the statistical analysis \citep{goldenberg}. 
Stochastic block models (SBMs), which originally stem from the social sciences \citep{white}, allow for a model-based clustering of nodes with similar connectivity behaviour, and hence are not only suitable for community detection, but for discovering any kind of connectivity pattern, see \citet{abbe, matias2014} for reviews of SBMs. 
Identification of the parameters in binary and weighted static SBMs was studied in \citet{Allman2009, Allman2011}. Recent contributions to  parameter estimation with maximum likelihood, variational maximum likelihood and Bayesian methods include \citet{celisse, bickel2013, Mariadassou15}.  

Instead of analyzing static snapshots of networks, some recent literature has been devoted to modeling network data which evolves through time \citep{holme}. 

In particular, \citet{Matias2016} propose a variant of the SBM in a discrete time framework, in which node states evolve with a Markovian dependence structure. They discuss identifiability for parametric models in the weighted case and introduce a variational EM-type algorithm for parameter estimation and a-posteriori clustering.  

In this note, complementary to the parametric identification for the dynamic SBM in \citet{Matias2016} we provide nonparametric identification results  in case of binary, finitely weighted and general edge states. We formulate conditions on the true parameters which guarantee actual point identification instead of mere generic identification. 
In particular, our results justify in terms of identification the applications in \citet{Matias2016} and in \citet{miele} to finitely weighted edges with $\kappa=3$ edge states. Our results also extend to inhomogeneous Markov chains. 	

A formal description of the model and the theoretical results are presented in Section \ref{sec:theory}. In Theorem \ref{th:nonparident} we consider nonparametric classes of distributions for weighted edges under the assumption that the true underlying distributions are linearly independent. Next, in Theorem \ref{th:smallstatic} we provide a result in the static setting for binary and finitely weighted edges with few edge states along the lines of the arguments in \citet{Allman2011}, but formulated and proved without referring to generic conditions. This is extended to the dynamic SBM in Theorem \ref{th:dynfewedgestates}. As corollary, we provide a convenient set of assumptions for identification in the dynamic SBM in the binary case, which we compare to the result in \citet{Matias2016a}. 
Remarks \ref{rem:inhomogeneous} and \ref{rem:inhomogeneous2} deal with extensions to inhomogeneous, nonstationary Markov chains. 
Section \ref{sec:sims} contains numerical illustrations for finitely weighted edges under scenarios to which Theorems \ref{th:nonparident} and \ref{th:dynfewedgestates} apply. Estimation is based on the variational EM algorithm from \citet{Matias2016}, for which we make the necessary modifications for the inhomogeneous case. 
Section \ref{sec:conclude} concludes, while proofs are deferred to Section \ref{sec:proofs}. 
 
We shall denote by $M \otimes N$ the Kronecker or tensor product of the matrices $M$ and $N$, by $M^{\otimes n}$, $n \in \N$, the $n$-fold Kronecker product of the matrix $M$ with itself, and by $M'$ the transpose of the matrix $M$. 
The entropy of a distribution $\mathbb{P}$ is denoted by $\mathcal{H}(\mathbb{P})$.

\section{The dynamic stochastic block model}\label{sec:theory}

Let $G=(V,E)$ be a complete, undirected graph without self-loops, where $V=\{1, \ldots,n\}$ is the finite set of nodes of cardinality $n \in \N$ and $E= \big\{ \{i,j\}, 1 \leq i \not= j \leq n\big\}$ is the set of undirected edges.  

The dynamic random graph mixture model or dynamic stochastic block model introduced in \cite{Matias2016} consists of latent random vectors $\boldsymbol{Z}^t=(Z_1^t,\ldots,Z_n^t)$ for the time points $t \in \boldsymbol{T} =\{1,\ldots,T\}$, $T \in \N$, where for each $i \in V$,  $(Z_i^t)_{1 \leq t \leq T}$ is a Markov chain with state space $\{1,\ldots, Q\}$, representing the node groups. The processes $(Z_i^t)_{1 \leq t \leq T}$ are assumed to be independent across individuals $i$, while their distributions are assumed not to depend on $i$. 

In the results below we suppose that the $(Z_i^t)_{1 \leq t \leq T}$ are time-homogeneous with ergodic transition matrix $\boldsymbol{\rho} = (\rho_{ql})_{q,l=1, \ldots, Q}$ and stationary initial distribution $\fpi = (\pi_q)_{q=1,\ldots, Q}$. However, the results extend to non-stationary, non-ergodic and even inhomogeneous Markov chains. 
Further, for each $t \in \boldsymbol{T}$ there is an observed random vector $\boldsymbol{X}^t=(X_e^t)_{e\in E}$, indexed by the edge set $E$ such that conditional on $\boldsymbol{Z}^1, \ldots, \boldsymbol{Z}^T$, the vectors $\boldsymbol{X}^1, \ldots, \boldsymbol{X}^T$ are independent, and conditional on $\boldsymbol{Z}^t$, the random variables $X_{\{i,j\}}^t$, $\{i,j\} \in E$ are independent, with the conditional distribution of $X_{\{i,j\}}^t$ depending only on $Z_i^t$ and $Z_j^t$. Further, the conditional distribution of $X_{\{i,j\}}^t$ given $Z_i^t$ and $Z_j^t$ does not depend on the edge $\{i,j\}$, it may depend on the time point $t$, though. We shall write $X_{\{i,j\}}^t = X_{ij}^t$, so that $X_{ji}^t= X_{ij}^t$, and call $X_{ij}^t$ the edge variables, and $Z_i^t$ the node states. 

Note that $(\boldsymbol{Z}^t,\boldsymbol{X}^t)_{t=1\ldots, T}$ is a hidden Markov model with state space $\{1, \ldots, Q \}^n$, transition probability matrix $\f{\rho}^{\otimes n}$ and stationary initial distribution $\f{\pi}^{\otimes n}$, and for each edge $\{i,j\}$,  $\big(({Z}_i^t,Z_j^t), X_{ij}^t\big)_{t=1\ldots, T}$ is a hidden Markov model with state space $\{1, \ldots, Q \}^2$, transition probability matrix $\f{\rho}^{\otimes 2}$ and stationary initial distribution $\f{\pi}^{\otimes 2}$. 

Suppose that the $X_{ij}^t$ take values in some measurable set $\mathcal{X}\subseteq \mathbb{R}^s$ containing $0$, called the set of edge states, such that the conditional distribution of $X_{ij}^t$ given $Z_i^t$ and $Z_j^t$ is given by
\begin{equation}\label{eq:conddistr}
\mu_{ql}^t:=\mathbb{P}(X_{ij}^t\cond Z_i^t=q,Z^t_j=l)=(1 - p_{ql}^t)\delta_0 + p_{ql}^t\, P_{ql}^t, \qquad  1\leq q, l \leq Q,
\end{equation}
where $\delta_0$ is the Dirac distribution at $0$, the $p_{ql}^t \in [0,1]$ are sparsity parameters, and the $P_{ql}^t$ are distributions on $\cX$ contained in a family $\cP$ of distributions which put measure $0$ for $\{0\}$, and we assume that $p_{ql}^t = p_{lq}^t$ and $P_{ql}^t = P_{lq}^t$. The value $0$ will be interpreted as the absence of an edge, while other values $x \in \cX$ correspond to the existence of a weighted edge. 
This general model is denoted by $\mathcal{M}(n,Q,T, \mathcal{X}, \cP)$. It has parameters $\boldsymbol{\rho}$, $p_{ql}^t$ and $P_{ql}^t$, $1 \leq q,l \leq Q$ and $1 \leq t \leq T$. 

For a {\sl binary edge states}, that is if $\cX = \{0,1\}$, we take $\cP = \{ \delta_1\}$, then only the parameters $\boldsymbol{\rho}$ and $p_{ql}^t$ are required, so that we may write $\mathcal{M}(n,Q,T, \{0,1\})$. 

For a general {\sl finite set of edges states}, that is, $\cX = \{0,1, \ldots, \kappa - 1\}$ with $\kappa \geq 2$, we use a parametrization which slightly differs from (\ref{eq:conddistr}), and parametrize the model in terms of the transition matrix $\boldsymbol{\rho}$ and of the full probability vectors of conditional edge states $\bp_{ql}^t(x) = \mathbb{P}(X_{ij}^t=x \cond Z_i^t = q, Z_j^t = l)$, $1 \leq q,l \leq Q$, $\{i,j\}\in E$, $x=0, \ldots, \kappa-1$, again assuming $\bp_{ql}^t(x) = \bp_{lq}^t(x)$. This corresponds to (\ref{eq:conddistr}) with sparsity parameter $p_{ql}^t = 1 - \bp_{ql}^t(0)$ and $P_{ql}^t(x) = \bp_{ql}^t(x)/(1-\bp_{ql}^t(0))$, $x \in \{1, \ldots, \kappa - 1\}$, the conditional distributions given $x \not=0$. We shall denote this model by $\mathcal{M}(n,Q,T, \{0,\ldots, \kappa -1\})$.  

If there is only a single time point, we obtain the static stochastic block model in \cite{Allman2011}, which is denoted by $\mathcal{M}(n,Q, \mathcal{X}, \cP)$, or $\mathcal{M}(n,Q, \{0,1\})$ resp.~$\mathcal{M}(n,Q, \{0,\ldots, \kappa -1\} )$ for binary resp.~finite edge set. In this model, the transition parameter $\boldsymbol{\rho}$ reduces to the distribution $\boldsymbol{\pi}$ of the node states at a single time point. 
 Apparently, the marginal distribution of  $(\boldsymbol{Z}^t,\boldsymbol{X}^t)$, $t \in \boldsymbol{T}$,  is contained in $\mathcal{M}(n,Q, \mathcal{X}, \cP)$ with parameters $\boldsymbol{\pi}$, $p_{ql}^t$ and  $P_{ql}^t$.

\section{Identification}

The parameterization in latent state models is subject to label swapping. Indeed, the distribution of the observables $(\boldsymbol{X}^t)_{t=1,\ldots, T}$ in the model $\mathcal{M}(n,Q,T, \mathcal{X}, \cP)$ is the same under the parameters $\rho_{ql}$, $p_{ql}^t$ and $P_{ql}^t$, as well as $\rho_{\sigma(q) \sigma(l)}$, $p_{\sigma(q) \sigma(l)}^t$ and $P_{\sigma(q) \sigma(l)}^t$, $1 \leq q,l \leq Q$ and $1 \leq t \leq T$, for any permutation $\sigma$ of the set of node states $\{1, \ldots, Q\}$. Note that the permutation $\sigma$ is {\sl global} in the sense that it must not depend on the time point $t$. 

Our aim is to discuss identification up to {\sl global label swapping}.

\subsection{Nonparametric models}\label{sec:nonpar}

We start with identification of the parameters in general nonparametric models with weighted edges. 

\begin{theorem}
	\label{th:nonparident}
	The parameters in the dynamic stochastic block model $\mathcal{M}(n,Q,T, \mathcal{X}, \cP)$ are identified from the distribution of $\boldsymbol{X}^1, \ldots, \boldsymbol{X}^T$ up to global label swapping if $n \geq 9$ and $T \geq 2$ provided that
	\begin{itemize}
		\item[a.] for every $t=1,\ldots,T$, the distributions $\mu_{ql}^t$ are linearly independent, $1 \leq q \leq l \leq Q$, and $p_{ql}^t>0$, 
		\item[b.] the entries of $\f{\pi}$ are pairwise distinct, or for $q=1,\ldots,Q-1$ the measures  $P^t_{qq}$ are stable over time, i.e.~$P_{qq}^t=P_{qq} \quad \forall \ t\in \boldsymbol{T}$.
	\end{itemize}
\end{theorem}
	An outline for the argument of the {\sl proof of Theorem \ref{th:nonparident}} is as follows. By ergodicity of $\boldsymbol{\rho} $, the entries of $\boldsymbol{\pi}$ are strictly positive. By assumption a.~and since $n \geq 9$, from \citet[theorem 15]{Allman2011} the marginal parameters $\boldsymbol{\pi}$, $p_{ql}^t$ and  $P_{ql}^t$ are identified from the distribution of $\boldsymbol{X}^t$ up to permutations $\sigma_t$ of the node states $\{ 1, \ldots, Q\}$ which, however, depend on $t$. Assumption $b.$~in the theorem then allows to align the node states globally, i.e.~to pass to a global permutation $\sigma$. Here, in the second case one observes that the distributions $P_{qq}$, $q=1, \ldots, Q$ are all distinct by assumption a.~in the theorem. It then remains to identify the transition matrix $\boldsymbol{\rho}$, which is the main part of the proof. 

\begin{remark}[{\sl Size of the set of edge states $\cX$}]
	First note that Theorem \ref{th:nonparident} can also be applied in case of a finite set of edges states, that is, model $\mathcal{M}(n,Q,T, \{0,\ldots, \kappa -1\})$. In terms of the $\bp_{ql}^t$, we require in a.~that $\bp_{ql}^t(0)<1$ and that $\bp_{ql}^t$, $1 \leq g \leq l \leq Q$ are linearly independent for fixed time point $t$. Note that this is only possible if the cardinality of the edge states $\kappa$ is at least as large as the number of random vectors, that is $\kappa \geq {Q+1 \choose 2}$, see \citet[theorem 14]{Allman2011}. In the second part of assumption b.~we require that the conditional probabilities $\bp_{qq}^t(x)/(1-\bp_{qq}^t(0))$, $x \in \{1, \ldots, \kappa - 1\}$, $q=1, \ldots, Q$, do not depend on $t$. \hfill $\diamond$
\end{remark}

\begin{remark}[{\sl Identification for nonstationary, inhomogeneous Markov chains}]\label{rem:inhomogeneous}
	Our argument even allows identification in case of inhomogeneous transitions of node states. Suppose that $\f{Z}_i$ is a possibly inhomogeneous Markov chain with some initial distribution $\fpi = (\pi_q)_{q=1,\ldots, Q}$ satisfying $\pi_q>0$, $q=1,\ldots, Q$, and irreducible transition matrices $\boldsymbol{\rho}^t = (\rho_{ql}^t)_{q,l=1, \ldots, Q}$, $t=2, \ldots, T$, giving the probability of transitions from time $t-1$ to $t$. Under the assumptions $a.$~and the second part of $b.$, that is, for $q=1,\ldots,Q-1$ the measures  $P^t_{qq} = P_{qq}$ are stable over time, the parameters $\fpi$, $p_{ql}^t$ and  $P_{ql}^t$, $t=1, \ldots, T$ and $\boldsymbol{\rho}^t$, $t=2, \ldots, T$, are still identified, see Section \ref{sec:proofs} for the argument. \hfill $\diamond$
\end{remark}

\begin{remark}[{\sl Assumptions only for the true parameters}]\label{rem:trueparameters}
	The assumptions $a.$~and $b.$ need only be placed on the true underlying parameters, identification is then achieved within the full class $\mathcal{M}(n,Q,T, \mathcal{X}, \cP)$. Indeed, the theorem relies on \citet[theorem 15]{Allman2011} which in turn is based on \citet[theorem 4a]{Kruskal} on the uniqueness of the factors in three-way arrays, see also \citet[theorem 3]{Rhodes2010}. These only require that the factors in one of two decompositions of a three-way array under comparison satisfy certain row restrictions. See \citet{Holzmann2016} for similar arguments in the case of nonparametric identification of hidden Markov models. \hfill $\diamond$
\end{remark}

\begin{remark}[{\sl Identification of the order $Q$}]
Consider the set of parameters on weighted dynamic SBMs with edges taking values in $\mathcal{X}$ and edge distributions $P_{ql}^t$ contained in some family $\cP$. If we let the order $Q$ vary, but only consider parameter constellations such that for given $Q$, the distributions $\mu_{ql}^t$ in (\ref{eq:conddistr}) satisfy assumption a.~of the theorem, that is, are linearly independent, then we claim that the parameters in this class, in particular the order $Q$ of the latent state space, are also identified. Indeed, an order $Q$-model can be regarded as having higher order e.g.~by choosing $p_{l,(q+1)} = p_{l,q}$, $P_{l,(q+1)} = P_{l,q}$, $l=1, \ldots, q$, and $p_{(q+1),(q+1)} = p_{q,q}$, $P_{(q+1),(q+1)} = P_{q,q}$, and the state $q$ in the Markov chain is randomly split into $q$ and $q+1$, from both of which the transitions have equal probability. However, if this is compared with a model on $Q+1$-states which does satisfy assumption a., then from Remark \ref{rem:trueparameters} it follows that the resulting distributions of observables are distinct.  
A similar remark applies in the situations of Theorems \ref{th:smallstatic} and \ref{th:dynfewedgestates}. \hfill $\diamond$
\end{remark}

\subsection{Binary and finitely weighted model}\label{sec:binary}

Next we consider identification in binary and finitely weighted models with few edge states. Assumption $a.$~in Theorem \ref{th:nonparident} cannot be satisfied in these cases. We develop an alternative set of conditions for point identification, complementing the generic identification results in \citet{dynSBM} for the dynamic and in \citet{Allman2011} for the static case.

Let us introduce a relevant object for identification. Given $m \in \N$ nodes  let $\boldsymbol{Z}=(Z_1,\ldots,Z_m)$ and $\boldsymbol{X}=(X_e)_{e \in E}$ follow the static stochastic block model $\mathcal{M}(m,Q, \{0,\ldots, \kappa -1\})$  and write $\bX = \bX_m$ to be explicit about the number $m$ of nodes. Consider the matrix $\cC_{m,Q,\kappa}$ consisting of the conditional distribution of $\boldsymbol{X}_m$ given $\boldsymbol{Z}$, so that 
\begin{align}\label{eq:thecondmatrix1}
\begin{split}
 \cC_{m,Q,\kappa}(\bz,\bx) & = \mathbb{P}(\bX_m=\bx \cond \bZ = \bz) = \prod_{1 \leq i < j \leq m} \bp_{q_i q_j}(x_{ij}), \\ 
 & \bx = (x_{ij})_{\{i,j\} \in E} \in \{0, \ldots, \kappa-1\}^E,\ \bz \in \{1, \ldots, Q\}^m,
 \end{split}
 \end{align}
and $\cC_{m,Q,\kappa}$ is of dimension $Q^m \times \kappa^{m \choose 2}$.

\medskip
First we state a point identification result for the static model under a full row-rank assumption on $\cC_{m,Q,\kappa}$. The proof is along the lines of \citet[theorem 14]{Allman2011}. Since the proof of that result is only sketched in \citet{Allman2011}, for completeness we provide a detailed proof of the following theorem.

\begin{theorem}\label{th:smallstatic}
Suppose that $m \geq 3$, and that $ Q, \kappa$,  and the conditional probabilities $\bp_{ql}(x)$, $1 \leq q \leq l \leq Q$ and $x \in \{0, \ldots, \kappa-1\}$ are such that the matrix $\cC_{m,Q,\kappa}$ has full row rank. Further assume that $\pi_q>0$, $q=1, \ldots, Q$.  Then $\boldsymbol{\pi}$ and $\bp_{ql}(x)$ are identified up to label swapping from the distribution of $\bX_n$ over the full model class $\mathcal{M}(n,Q, \{0, \ldots, \kappa-1\})$ provided that $n \geq m^2$.
\end{theorem}

\begin{remark}[{\sl Assumptions only for the true parameters}]
	The assumption of a full row rank of $\cC_{m,Q,\kappa}$ need only hold for the true conditional probabilities $\bp_{ql}(x)$, identification is then achieved over the full class $\mathcal{M}(n,Q, \{0, \ldots, \kappa-1\})$, similarly for the restriction $\pi_q>0$.  \hfill $\diamond$
\end{remark}

Theorem \ref{th:smallstatic} provides conditions for identification in case of a small set of edge states that can be verified by a rank computation with symbolic algebra software.
Using MATLAB, we construct the symbolic matrix $\cC_{m,Q,\kappa}$ for fixed $Q$ and $\kappa\leq 10$ in case that all free parameters $\bp_{ql}(x)$, $1 \leq q \leq l \leq Q$, $x \in \{0, \ldots, \kappa-2\}$ are distinct.
Iteratively, we determine for several combinations of $\kappa$ and $Q$ the minimal number $m\in \mathbb{N}$ such that the matrix $\cC_{m,Q,\kappa}$ has full row-rank. Thus, the model on $n\geq m^2$ nodes is identifiable if additionally $\pi_q>0$, $q=1, \ldots, Q$. Table \ref{Tab:Matlab} contains the results for $2 \leq Q \leq 5$. 
For example, from Table \ref{Tab:Matlab} it can be inferred that a binary model ($\kappa=2$) with $Q=3$ latent node states is identifiable as soon as the network consists of at least $m^2=5^2=25$ nodes.
Additional to the results given by Table \ref{Tab:Matlab}, we showed that also for all combinations $\kappa \geq Q$ with $6\leq Q \leq 10$, a minimum of 3 nodes is enough to ensure a matrix of full row-rank.
\begin{table}
\centering 
\setlength{\tabcolsep}{5mm}
\begin{threeparttable}
\caption{The minimal number $m$ for combinations of $\kappa$ and $Q$ such that the matrix $\cC_{m,Q,\kappa}$ has full row-rank for every set of probability vectors $\bp_{ql}$, $1 \leq q \leq l \leq Q$ for which all of the free parameters $\bp_{ql}(x)$, $1 \leq q \leq l \leq Q$, $x \in \{0, \ldots, \kappa-2\}$ are distinct.}
\label{Tab:Matlab}
\begin{tabular}{|c|c|c|c|c|c|c|c|c|c|}
\hline 
\backslashbox{$Q$}{$\kappa$} & 2 & 3 & 4 & 5 & 6 & 7 & 8 & 9 & 10 \\ 
\hline 
2 & 4 & 3 & 3 & 3 & 3 & 3 & 3 & 3 & 3 \\ 
\hline 
3 & 5 & 3 & 3 & 3 & 3 & 3 & 3 & 3 & 3 \\ 
\hline 
4 & -- & 4 & 3 & 3 & 3 & 3 & 3 & 3 & 3 \\ 
\hline 
5 & -- & -- & 4 & 3 & 3 & 3 & 3 & 3 & 3 \\ 
\hline 
\end{tabular} 
\end{threeparttable}
\end{table}
\medskip

From Theorem \ref{th:smallstatic} we conclude the following general result. 

\begin{theorem}\label{th:dynfewedgestates}
	Consider the dynamic model $\mathcal{M}(n,Q,T, \{0,\ldots, \kappa-1\})$. Suppose that $m\geq 3$ and that $Q$, $\kappa$,  and the conditional probabilities $\bp_{ql}^t$, $t \in \boldsymbol{T}$ are such that the matrices $\cC_{m,Q,\kappa}^t$ formed as in (\ref{eq:thecondmatrix1}) with $\bp_{ql}^t$ have full row rank for each $t$, and that the number of nodes satisfies $n \geq m^2$. 
	Further suppose that the transition matrix $\f{\rho}$ is ergodic and has full rank, that $T \geq 3$, and that the distributions $\bp_{qq}^t = \bp_{qq}$ are stable in time and distinct for $q=1, \ldots, Q$.  
	
	Then the true parameters $\boldsymbol{\rho}$ and $\bp_{ql}^t$ are identified within the full model class \linebreak $\mathcal{M}(n,Q,T, \{0,\ldots, \kappa-1\})$ up to global label swapping. 
\end{theorem}

Again the conditions need only be placed on the true parameters, identification is then obtained within the full class $\mathcal{M}(n,Q,T, \{0,\ldots, \kappa-1\})$. As in Theorem \ref{th:nonparident} the proof proceeds via the marginal models $(\bZ^t, \bX^t)$ using Theorem \ref{th:smallstatic}. The identification of the transition matrix requires new arguments, however. 

Instead of assuming that the distributions $\bp_{qq}^t = \bp_{qq}$ are stable in time, one could impose other assumptions that allow for global assignment from local assignments of the node states, e.g.~that the entries of $\f{\pi}$ are pairwise distinct. 



\begin{remark}[{\sl Identification for nonstationary, inhomogeneous Markov chains}]\label{rem:inhomogeneous2}
	The argument can also be extended to include identification in case of inhomogeneous transitions of node states. Consider the assumptions in Theorem \ref{th:dynfewedgestates}, but instead of assuming a single ergodic transition matrix $\boldsymbol{\rho}$ and stationary marginal distribution $\boldsymbol{\pi}$, consider the inhomogeneous setting as in Remark \ref{rem:inhomogeneous}, and additionally assume that all transition matrices $\boldsymbol{\rho}^t$ have full rank. Then the parameters $\fpi$, $\bp_{ql}^t$, $t=1, \ldots, T$ and $\boldsymbol{\rho}^t$, $t=2, \ldots, T$, are still identified, see Section \ref{sec:proofs} for the argument.  \hfill $\diamond$
\end{remark}

Specifically, for the binary model we obtain the following result. 

\begin{cor}
Consider the dynamic model $\mathcal{M}(n,Q,T, \{0,1\})$ with binary edge state. Suppose that the transition matrix $\f{\rho}$ has full rank, and that $T \geq 3$. 
\begin{itemize}
	\item $Q=2$:\quad Assume that card$\{	p_{11}^t,p_{12}^t,p_{22}^t\} = 3$ for all $t \in \boldsymbol{T}$ and that $\rho_{jj} \in (0,1)$, $j=1,2$. If $\rho_{12}\not=\rho_{21}$ or if $p_{11}^t = p_{11}$ for all $t \in \boldsymbol{T}$ does not depend on $t$, then the parameters are identified within $\mathcal{M}(n,2,T, \{0,1\})$ provided that $n \geq 16$. 
	\item $Q=3$:\quad Assume that card$\{	p_{11}^t,p_{12}^t,p_{13}^t, p_{22}^t,p_{23}^t,p_{33}^t\} = 6$ for all $t \in \boldsymbol{T}$ and that $\f{\rho}$ is ergodic. If $p_{11}^t = p_{11}$, $p_{22}^t = p_{22}$ for all $t \in \boldsymbol{T}$ does not depend on $t$, then the parameters are identified within $\mathcal{M}(n,3,T, \{0,1\})$ provided that $n \geq 25$. 
\end{itemize}
\end{cor}

\citet{Matias2016a} also discuss the binary case. Their argument requires the generic condition that the average edge probabilities $\sum_{l=1} \pi_l \, p_{ql}^t$, $q=1, \ldots, Q$ are all distinct, which is not necessary for our result. 

Further, in the binary model, \citet[Theorem 2]{Allman2011} show that $\cC_{m,Q,2}$ generically has full row rank if $m$ is not too small with respect to $Q$. More precisely, if 
\begin{equation}
m \geq
\begin{cases}
     Q-1+\frac{1}{4}(Q+2)^2 & \text{if  $Q$ is even,} \\
     Q-1+\frac{1}{4} (Q+1)(Q+3) & \text{if $Q$ is odd. }
   \end{cases}   
\end{equation}
This is the only step involving a generic argument. Once it can be shown that the matrices $\cC^t_{m,Q,2}$ have full row rank, $t=1, \ldots, T$, the transition matrix is also identified.

\section{Estimation and numerical illustrations}\label{sec:sims}
In this section we illustrate numerically the estimation performance of the variational EM algorithm from \citet{Matias2016}, which is implemented in the \textit{dynsbm} package (available on CRAN). After clarifying identification it is of some interest to present numerical results in the weighted case which is particularly relevant in real-data applications, see \citet{Matias2016}, who however only investigate the binary situation in simulations.   

We consider two simulation settings in which the parameters are identified by Theorems \ref{th:nonparident} and \ref{th:dynfewedgestates}, respectively.
In both scenarios we make repeated simulations on $100$ networks of different sizes: A rather small network size of 150 nodes and few time points to show that estimation works even under minimal conditions required by Theorems \ref{th:nonparident} and \ref{th:dynfewedgestates} and a larger network size of 1000 nodes. Moreover, we show simulation results in case that the underlying Markov chain is inhomogeneous. Below we briefly discuss the necessary changes to the variational EM algorithm, which we implemented by modifying the code in the \textit{dynsbm} package. As recommended in \citet{Matias2016}, we choose $25$ starting points for the iterative algorithm for each network to ensure accurate estimation results.
All simulations were performed on a Unix workstation with 16 GB RAM and an eight-core Xeon E5-1620 v3 processor. 

\subsection{Estimation of inhomogeneous dynamic SBMs}

\citet{Matias2016} introduced a variational EM Algorithm for estimation of dynamic SBMs, which can be extended to inhomogeneous Markov chains.
The complete data log-likelihood of an inhomogeneous stochastic block model is of the form
\begin{small}
\begin{align*}
\mathbb{P}(\boldsymbol{X}_n,\boldsymbol{Z})&=\prod_{i=1}^n \prod_{q=1}^Q \pi^1_q \mathds{1}_{q}(Z_i^1) \prod_{t=2}^T \prod_{i=1}^n \prod_{1\leq q,l \leq Q} \mathds{1}_{q}(Z_i^{t-1})\mathds{1}_l(Z_i^t) \rho_{ql}^t \\
&\times \prod_{t=1}^T \prod_{1 \leq i,j \leq n} \prod_{1 \leq q,l \leq Q} \mathds{1}_{q}(Z_i^t) \mathds{1}_l(Z_j^t) \mu_{ql}^t(X_{ij}^t),
\end{align*}
\end{small}
with model parameters $\theta=(\fpi, (\boldsymbol{\rho}^t), (\mu_{ql}^t)).$\\\
The variational EM algorithm aims at optimizing the function
$$J(\theta, \mathbb{Q})= \mathbb{E}_{\mathbb{Q}}[\log\mathbb{P}(\f{X},\f{Z})] +  \mathbb{E}_{\mathbb{Q}}[\log\mathbb{Q}(\f{Z})	]$$ over $\mathbb{Q}\in \mathcal{D}$, where $\mathcal{D}$ is chosen as the class of factorizable distributions with the Markov property $$\mathcal{D}:=\left\{\mathbb{Q} ~ \cond~ \mathbb{Q}(\f{Z})= \prod_{i=1}^n \mathbb{Q}_i(Z^1_i)\prod_{t=1}^T \mathbb{Q}_i(Z^t_i\cond Z^{t-1}_i)\right\}.$$
The distributions $\mathbb{Q}_i$ are assumed to be further factorizable into 
\begin{align*}
\mathbb{Q}_i(Z_i^1)&=\prod_{q=1}^Q {(\lambda_q^i)}^{\mathds{1}_q(Z_i^1)},\qquad
\mathbb{Q}_i(Z^t_i \cond Z^{t-1}_i)= \prod_{1 \leq q,l \leq Q} { (^t\lambda_{q,l}^i)}^{\mathds{1}_q(Z_i^{t-1})\mathds{1}_l(Z_i^t)},	
\end{align*}
and are parameterized by the variational parameters ${^1\lambda_q^i}$ and ${ ^t\lambda_{q,l}^i}$ which correspond to 
\begin{align*}
{ ^1\lambda_q^i} &\approx \mathbb{P}(Z_i^1=q \cond \f{X}), \quad 
{ ^t\lambda_{q,l}^i} \approx \mathbb{P}(Z_i^t=l \cond Z_i^{t-1}=q, \f{X}) \quad \text{for } t \geq 2.
\end{align*}

The marginal components of a distribution $\mathbb{Q} \in \mathcal{D}$, $ ^t\delta_q^i \approx \mathbb{P}(Z_i^t=q\cond \f{X}),$ are computed recursively.

The VEM-Algorithm iteratively updates the variational parameters given the model parameters in a variational E-step, and the model parameters given the variational parameters in an M-step. Explicit formulas for the fixed point equation for the E-step as well as for the sparsity parameter in the M-step 	$p_{ql}^t$ can be found in Proposition 2 of \citet{Matias2016}.
Further, in the inhomogeneous case only the M-step for the transition matrices and the group proportions changes, so that distinct parameters for every time point are calculated as
\begin{align}
\pi_{q}^1&= \frac{1}{n} \sum_{i=1}^n {(^1\delta_q^i)}, \qquad 
\rho_{ql}^t \propto \sum_{i=1}^n { (^{t-1} \delta_q^i)} \cdot { (^t\lambda_{q,l}^i)}, \quad t= 2, \ldots,T
\end{align}
 instead of averaging over all time points.

\subsection{Scenario 1: Large set of edge states}
\begin{figure}
\subfigure[$n=150$]{\includegraphics[width=0.5\textwidth]{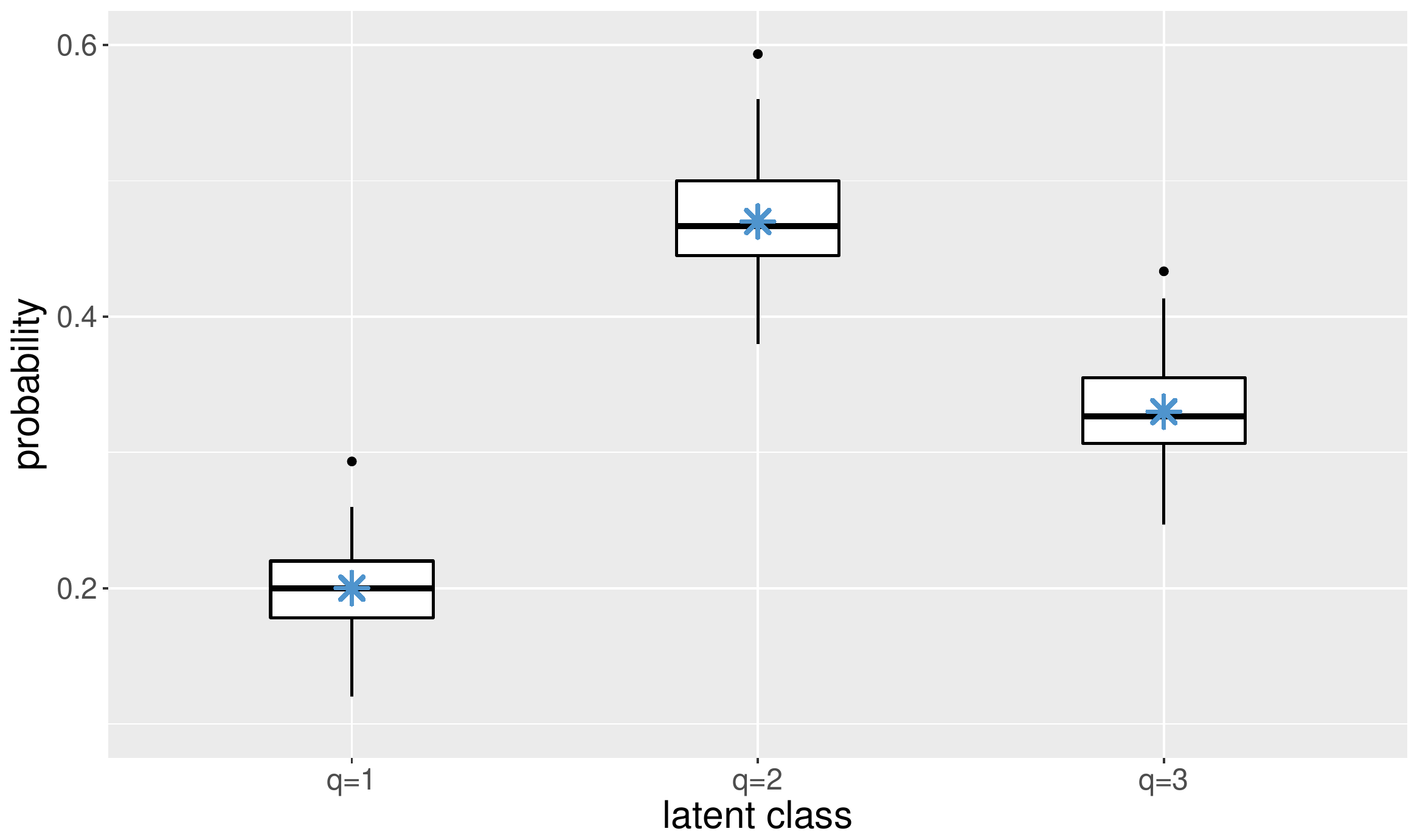}}\hfill
\subfigure[$n=1000$]{\includegraphics[width=0.5\textwidth]{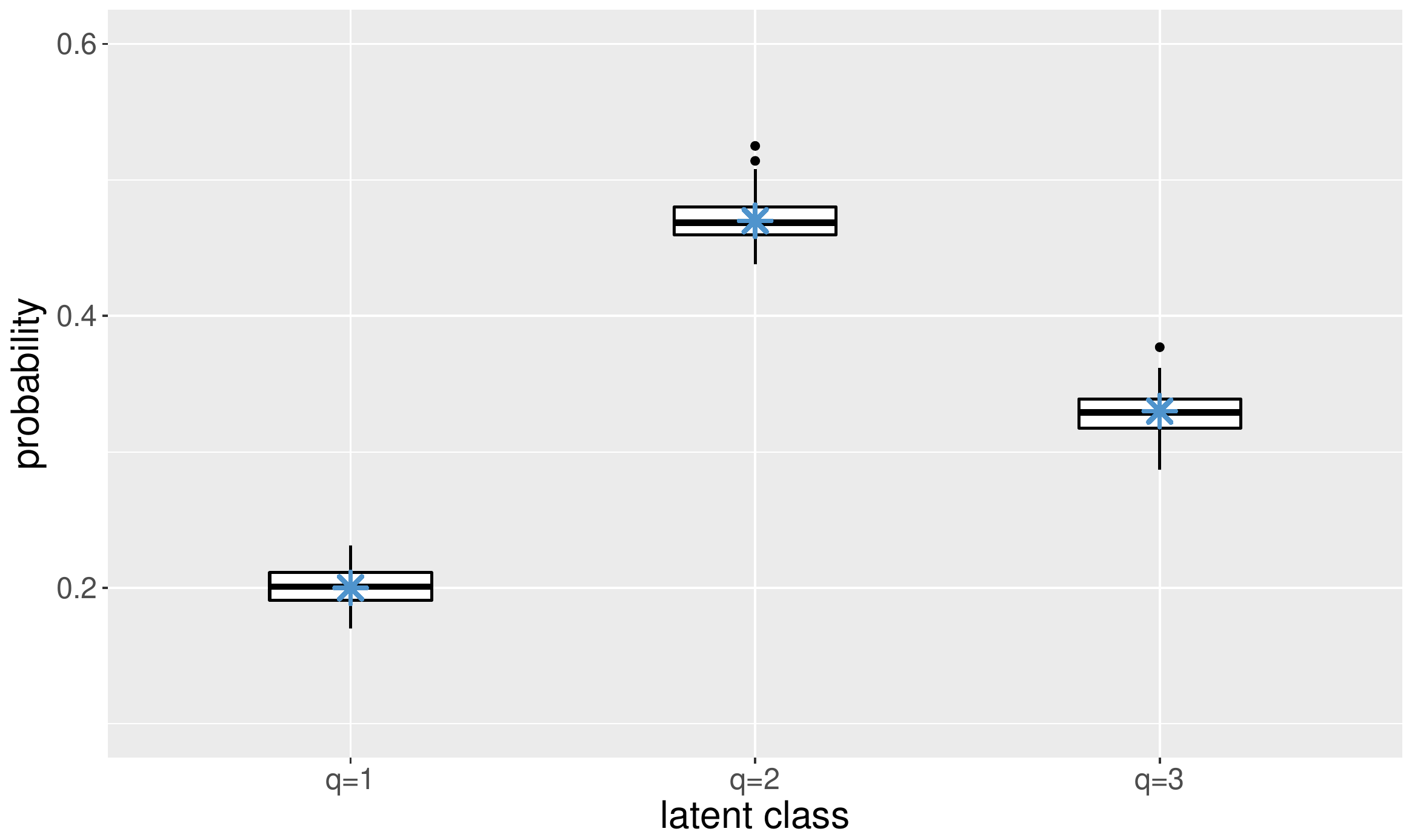}
}\hfill
\caption{Boxplots of estimation results of the starting group proportions $\fpi$ in Scenario 1 for 100 simulated networks. True underlying parameters are visualized as blue asterisks.  }
\label{fig:Sc1_pi}
\end{figure}

\begin{figure}
	\subfigure[$n=150$]{\includegraphics[width=0.5\textwidth]{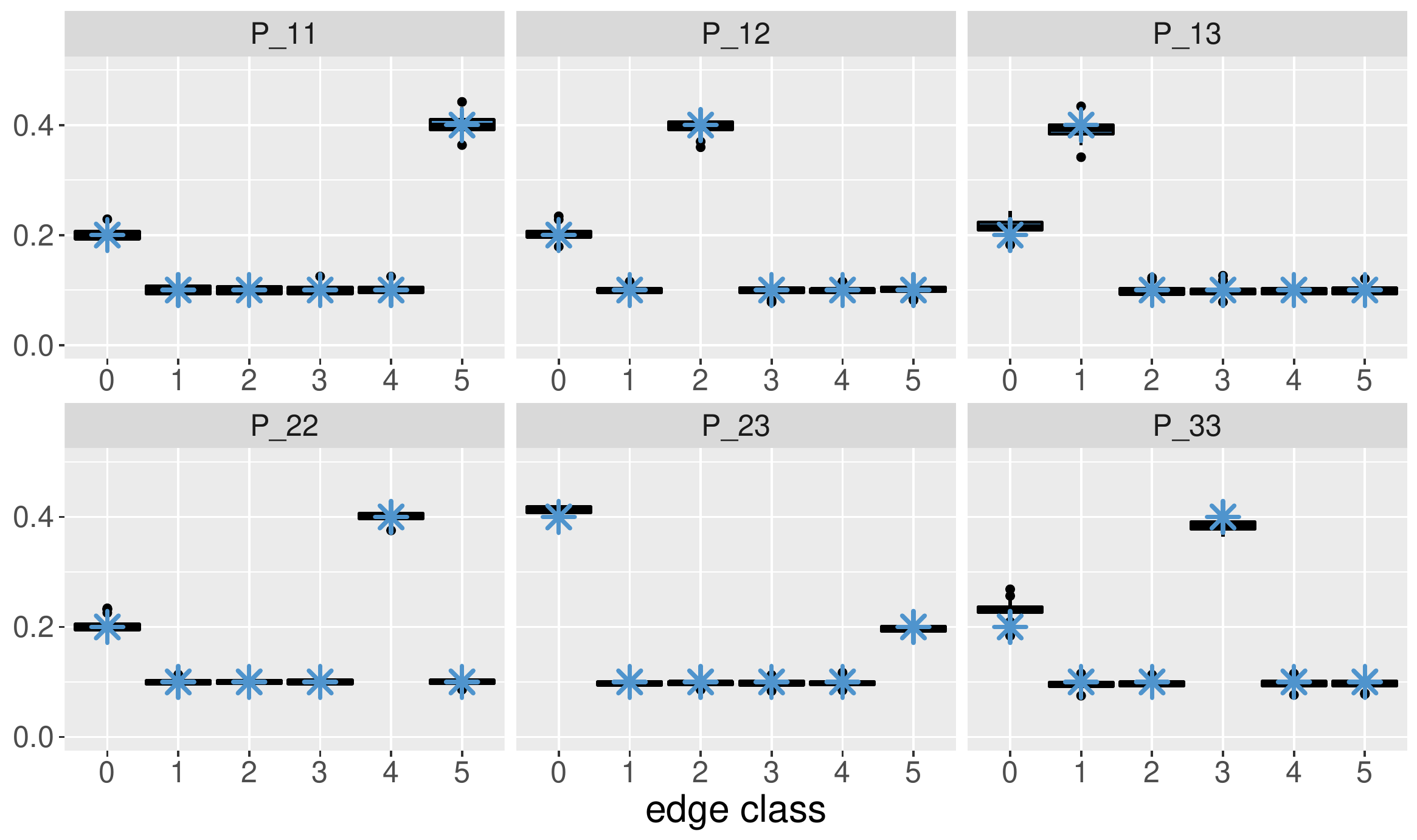}
	}\hfill
	\subfigure[$n=1000$]{\includegraphics[width=0.5\textwidth]{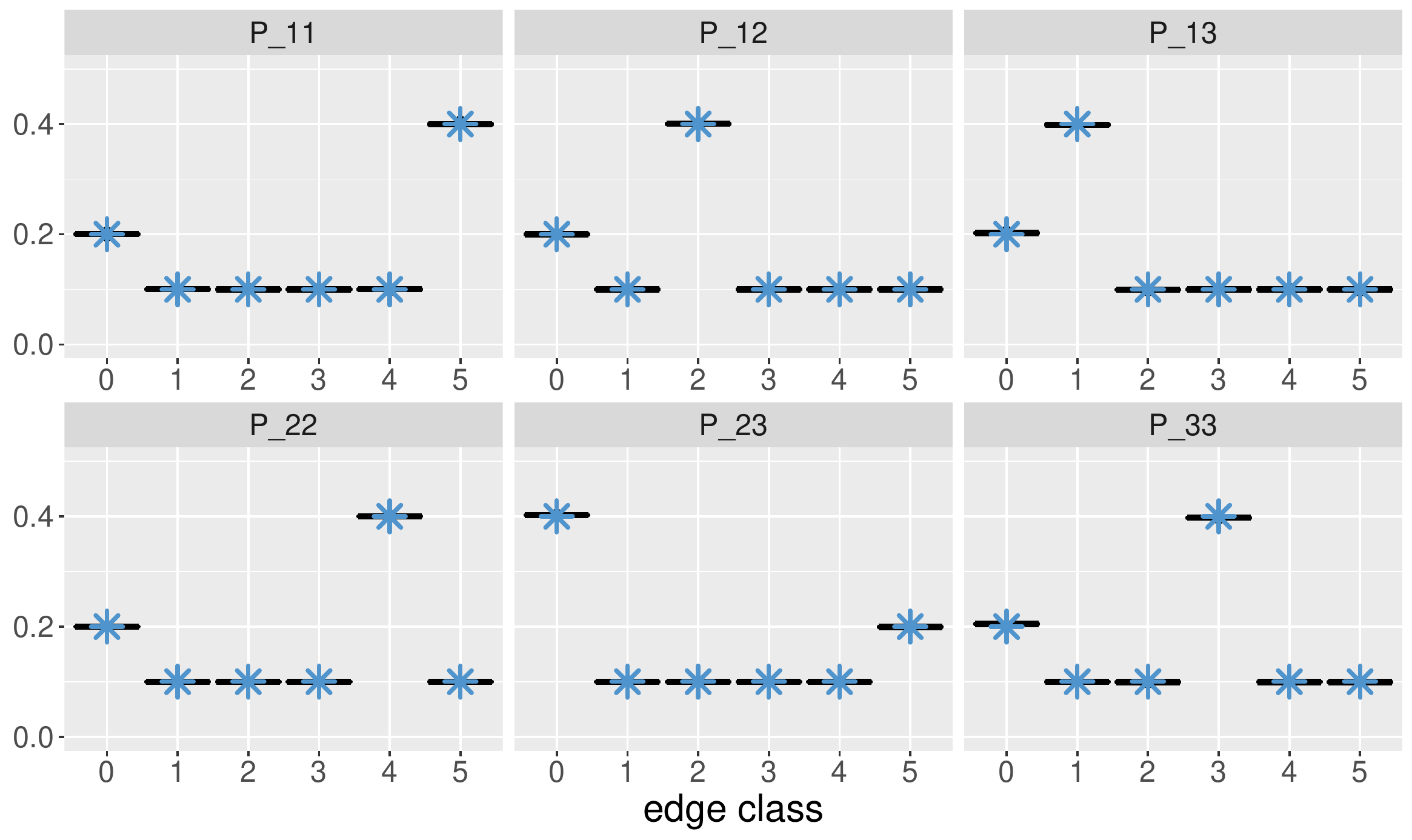}		
	}\hfill
\caption{Boxplots of estimation results of the distributions $\bp_{ql}$ in Scenario 1 for 100 simulated networks. True underlying parameters are visualized as blue asterisks. }
\label{fig:Sc1_Gamma}
\end{figure}

\begin{figure}
	\subfigure[$n=150$]{\includegraphics[width=0.5\textwidth]{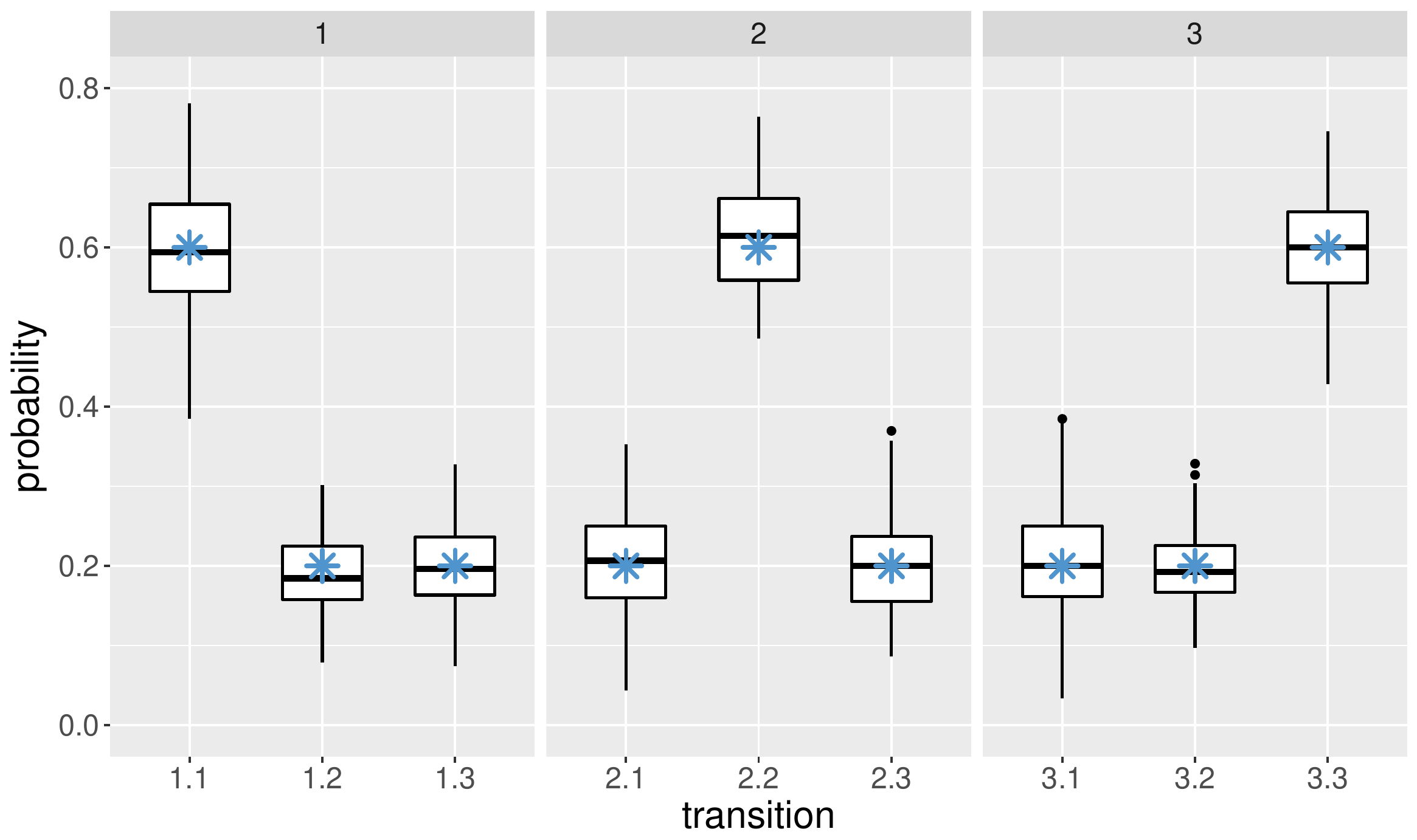}
	}\hfill
	\subfigure[$n=1000$]{\includegraphics[width=0.5\textwidth]{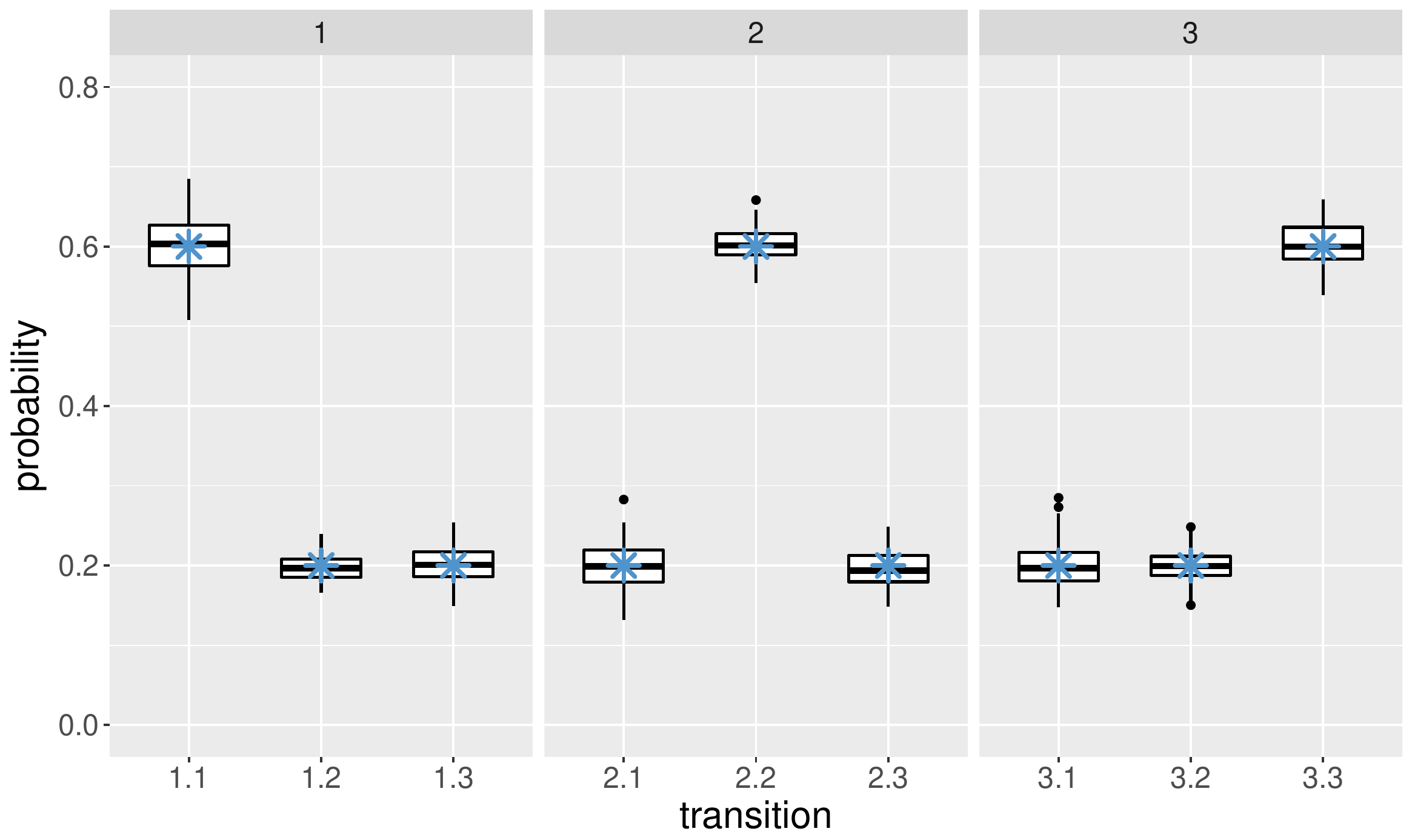}		
	}\hfill
\caption{Boxplots of estimation results of the transition matrix $\boldsymbol{\rho}$ in Scenario 1 for 100 simulated networks. True underlying parameters are visualized as blue asterisks.}
e2\label{fig:Sc1_Trans}
\end{figure}

In Scenario 1 we simulate random graphs on $n=150$ and on $n=1000$ nodes with $Q=3$ latent groups for $T=2$ different time points. To choose linearly independent distributions $\bp_{ql}$ and thus meet the conditions of Theorem \ref{th:nonparident} we consider $\kappa=\binom{4}{2}=6$ edge classes.
Global assignment of groups is ensured by taking all the edge distributions $\bp_{ql}$ stable over time, and we choose
{\small 
\begin{equation}\label{eq:tpmsim}
\begin{pmatrix}
\bp_{11} \\ 
\bp_{22} \\ 
\bp_{33} \\ 
\bp_{12} \\ 
\bp_{13} \\ 
\bp_{23}
\end{pmatrix} 
=
\begin{pmatrix}
0.2 & 0.1 & 0.1 & 0.1 & 0.1 & 0.4 \\ 
0.2 & 0.1 & 0.1 & 0.1 & 0.4 & 0.1 \\ 
0.2 & 0.1 & 0.1 & 0.4 & 0.1 & 0.1 \\ 
0.2 & 0.1 & 0.4 & 0.1 & 0.1 & 0.1 \\ 
0.2 & 0.4 & 0.1 & 0.1 & 0.1 & 0.1 \\ 
0.4 & 0.1 & 0.1 & 0.1 & 0.1 & 0.2
\end{pmatrix}, \quad \boldsymbol{\rho}=\begin{pmatrix}
0.6 & 0.2 & 0.2 \\ 
0.2 & 0.6 & 0.2 \\ 
0.2 & 0.2 & 0.6
\end{pmatrix}  
\end{equation}
}
and $\fpi=(0.2,0.33,0.47)$. Note that this is not the stationary distribution of $\boldsymbol{\rho}$, see Remark \ref{rem:inhomogeneous}.
In the repeated simulations, as the groups are only identified up to label swapping we need to determine an appropriate labeling in each run, which we base on the edge distributions $\bp_{ql}$.
The estimates of the parameters are visualized as boxplots in Figures \ref{fig:Sc1_pi}--\ref{fig:Sc1_Trans} for both sizes of the network, where the true parameter values are presented as blue asterisks. The edge distributions are estimated rather precisely, while there is more variability in the estimates of the transition probabilities and the initial group proportions. As expected, all estimates become more exact when the network size increases.  
The runtime is $\sim$10 minutes for $n=150$ and $\sim$18 hours for the larger network with $n=1000$.

\begin{figure}
	\subfigure[$\boldsymbol{\rho}^2$]{\includegraphics[width=0.5\textwidth]{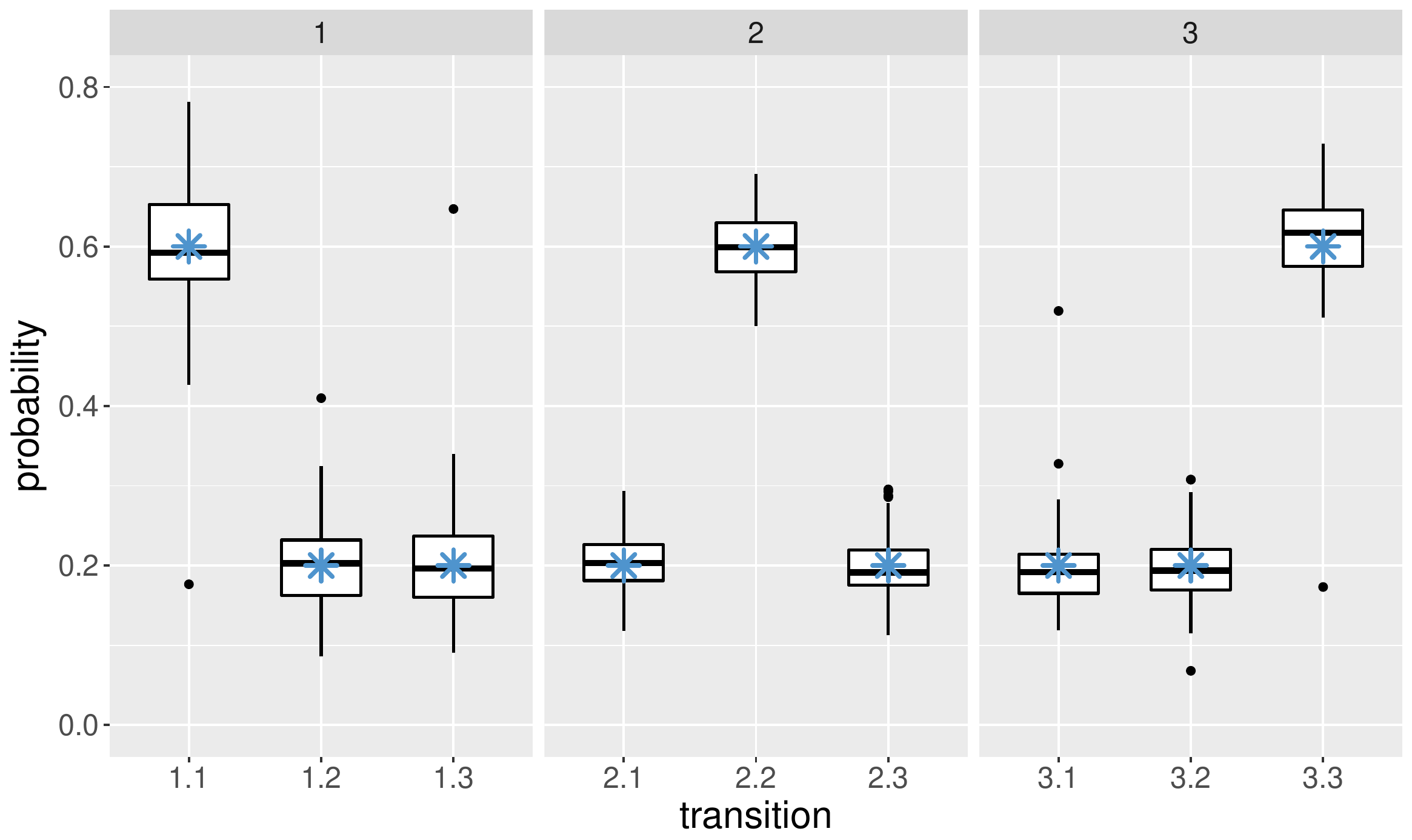}
	}\hfill
	\subfigure[$\boldsymbol{\rho}^4$]{\includegraphics[width=0.5\textwidth]{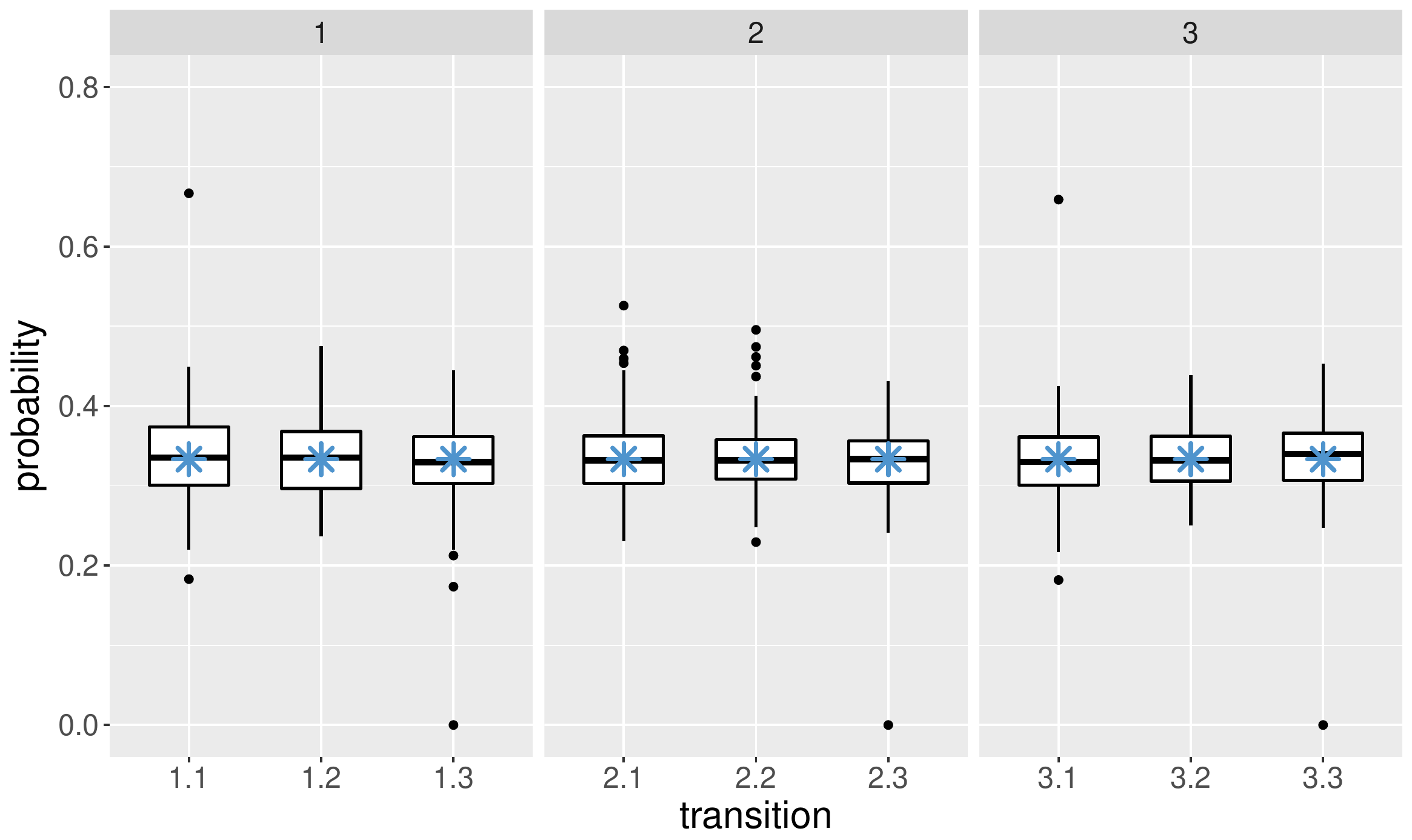}		
	}\hfill
\caption{Boxplots of estimation results of the transition matrices $\boldsymbol{\rho}^2$ and $\boldsymbol{\rho}^4$ in Scenario 1 for 100 simulated networks. True underlying parameters are visualized as blue asterisks.}
\label{fig:Sc2_Trans1_inhomo}
\end{figure}

Additionally we consider Scenario 1, that is $Q=3$, $\kappa = 6$ with edge distributions as in (\ref{eq:tpmsim}) and stable in time, but over $T=4$ time points and with an inhomogeneous Markov chain, where the transition matrices are given by
\begin{equation*}
\boldsymbol{\rho}^2=\begin{pmatrix}
0.6 & 0.2 & 0.2 \\ 
0.2 & 0.6 & 0.2 \\ 
0.2 & 0.2 & 0.6
\end{pmatrix},\quad  
\boldsymbol{\rho}^3=\begin{pmatrix}
0.15 & 0.15 & 0.7 \\ 
0.7 & 0.15 & 0.15 \\ 
0.15 & 0.7 & 0.15
\end{pmatrix},\quad 
\boldsymbol{\rho}^4=\begin{pmatrix}
\frac{1}{3} & \frac{1}{3} & \frac{1}{3} \\[6pt]
\frac{1}{3} & \frac{1}{3} & \frac{1}{3} \\[6pt] 
\frac{1}{3} & \frac{1}{3} & \frac{1}{3}
\end{pmatrix}.
\end{equation*}
We consider the medium size of $n=300$ nodes for the network. 
In Figure \ref{fig:Sc2_Trans1_inhomo} we display the estimation results for $\boldsymbol{\rho}^2$ and $\boldsymbol{\rho}^4$, for $\boldsymbol{\rho}^3$ these are similar. There is not much difference compared to the estimation quality in the homogeneous case. Here the runtime was $\sim$2 hours. 

\subsection{Scenario 2: Small set of edge states}
In Scenario 2 we simulate random graphs on $n=150$ and on $n=1000$ nodes with $Q=3$ latent groups for $T=3$ different time points and transition matrix $\boldsymbol \rho$ as in (\ref{eq:tpmsim}), which is of full rank. According to Table \ref{Tab:Matlab}, $\kappa=3$ different edge classes and pairwise distinct parameters $\bp_{ql}$ are enough to ensure identifiability based on Theorem  \ref{th:dynfewedgestates}. As starting group proportions we take $\boldsymbol{\pi}=\big(\frac{1}{3},\frac{1}{3},\frac{1}{3}\big)$, the stationary distribution of $\boldsymbol{\rho}$. Again, global assignment of groups is ensured by taking the edge distributions $\bp_{ql}$ stable over time, and we choose
$$\begin{pmatrix}
\bp_{11} \\ 
\bp_{22} \\ 
\bp_{33} \\ 
\bp_{12} \\ 
\bp_{13} \\ 
\bp_{23}
\end{pmatrix} 
=
\begin{pmatrix}
0.1 & 0.55 & 0.35 \\ 
0.2 & 0.45 & 0.35\\ 
0.3 & 0.35 & 0.35  \\ 
0.4 & 0.25 & 0.35  \\ 
0.5 & 0.15 & 0.35 \\ 
0.6 & 0.05 & 0.35 
\end{pmatrix} .
$$

\begin{figure}
	\subfigure[$n=150$]{\includegraphics[width=0.5\textwidth]{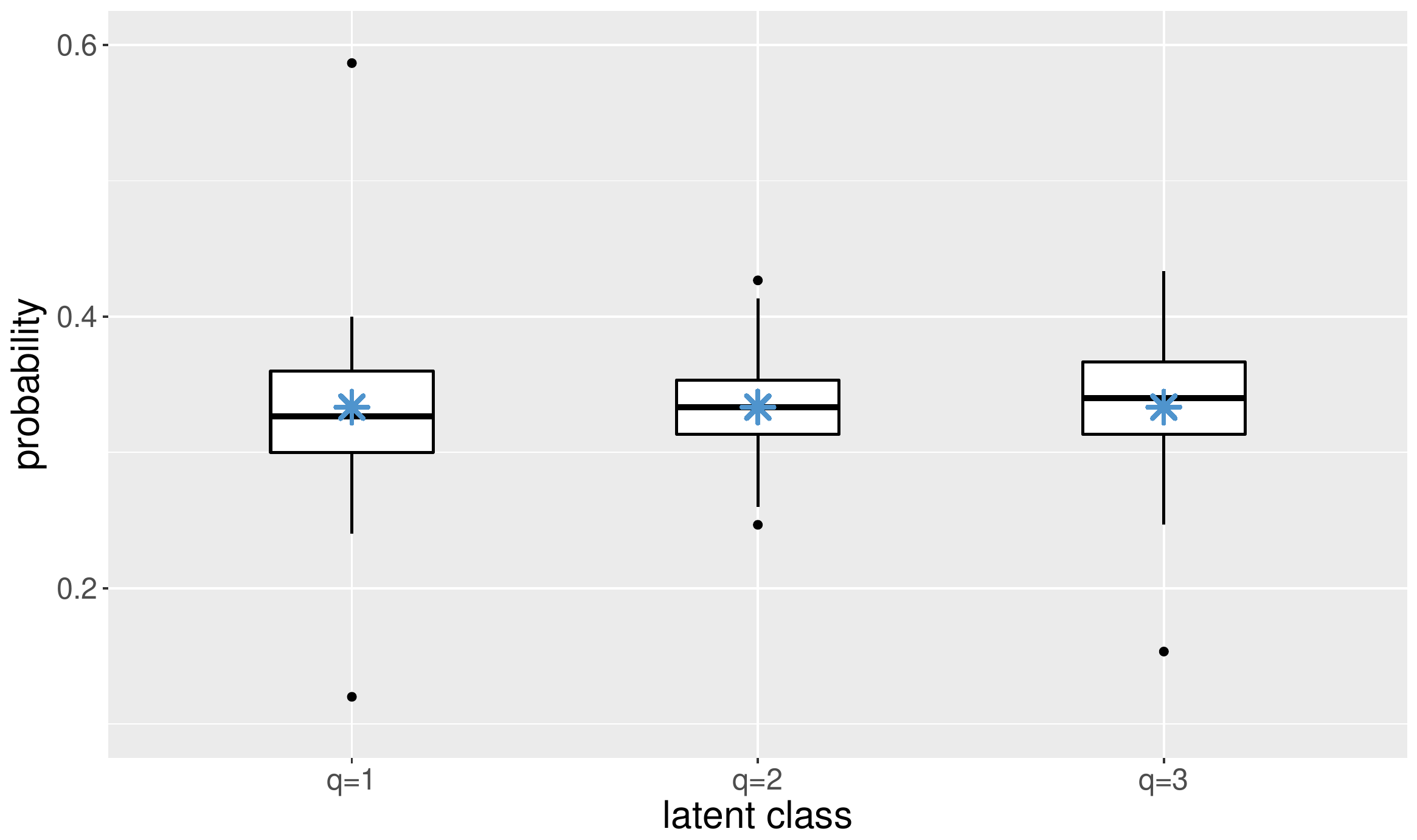}
	}\hfill
	\subfigure[$n=1000$]{\includegraphics[width=0.5\textwidth]{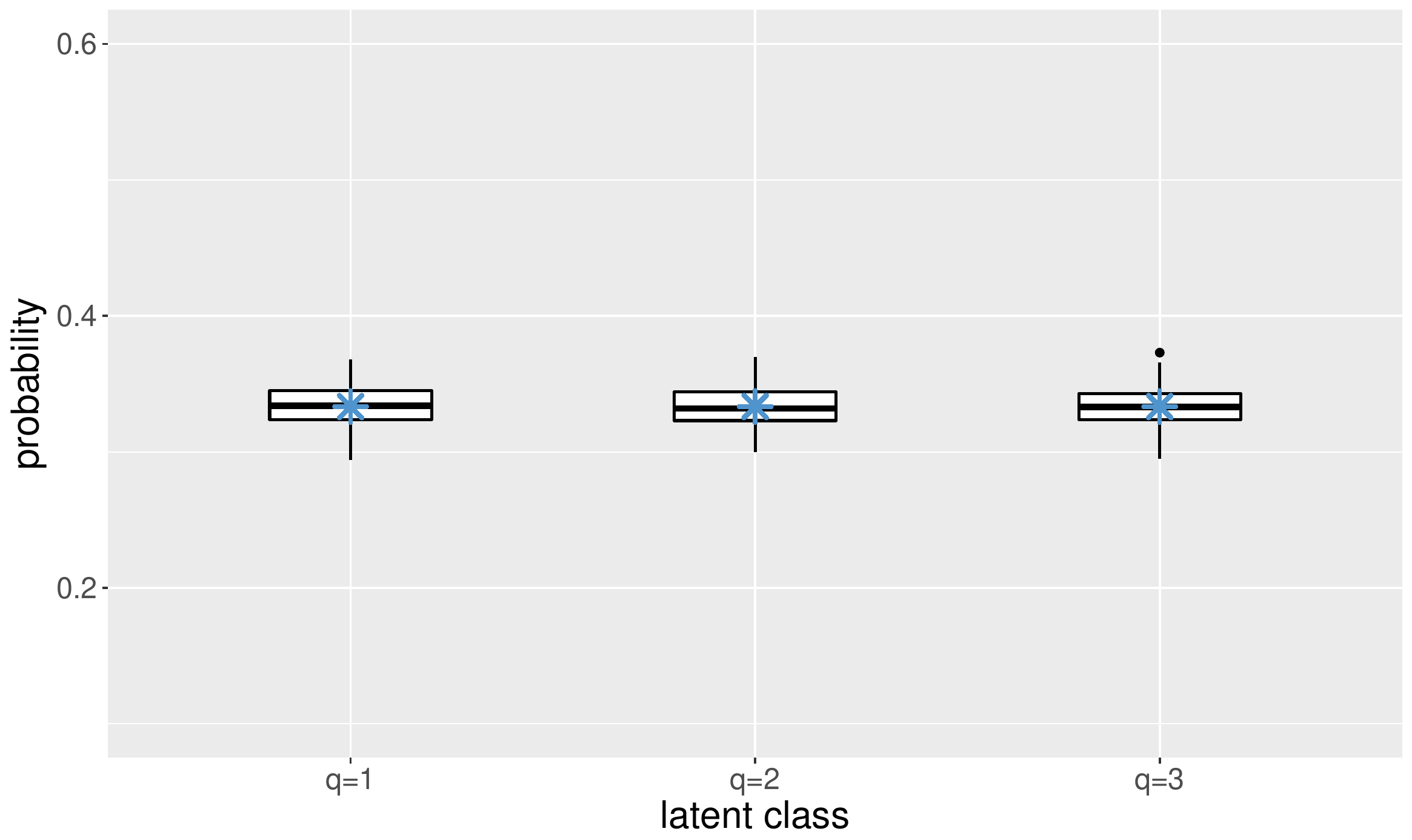}
	}\hfill
	\caption{Boxplots of estimation results of the starting group proportions $\fpi$ in Scenario 2 for 100 simulated networks. True underlying parameters are visualized as blue asterisks.}
	\label{fig:Sc2_pi}
\end{figure}

\begin{figure}
	\subfigure[$n=150$]{\includegraphics[width=0.5\textwidth]{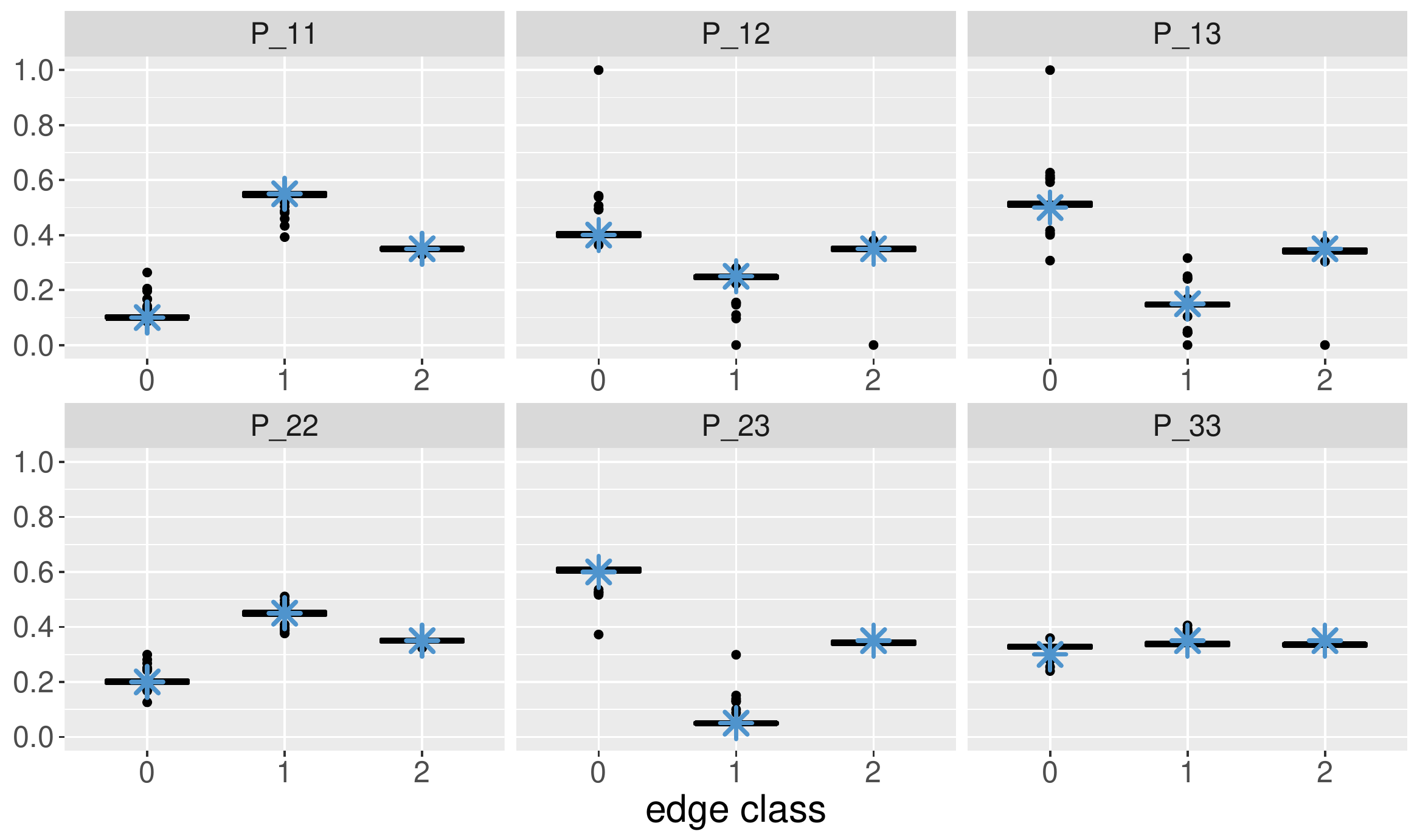}
	}\hfill
	\subfigure[$n=1000$]{\includegraphics[width=0.5\textwidth]{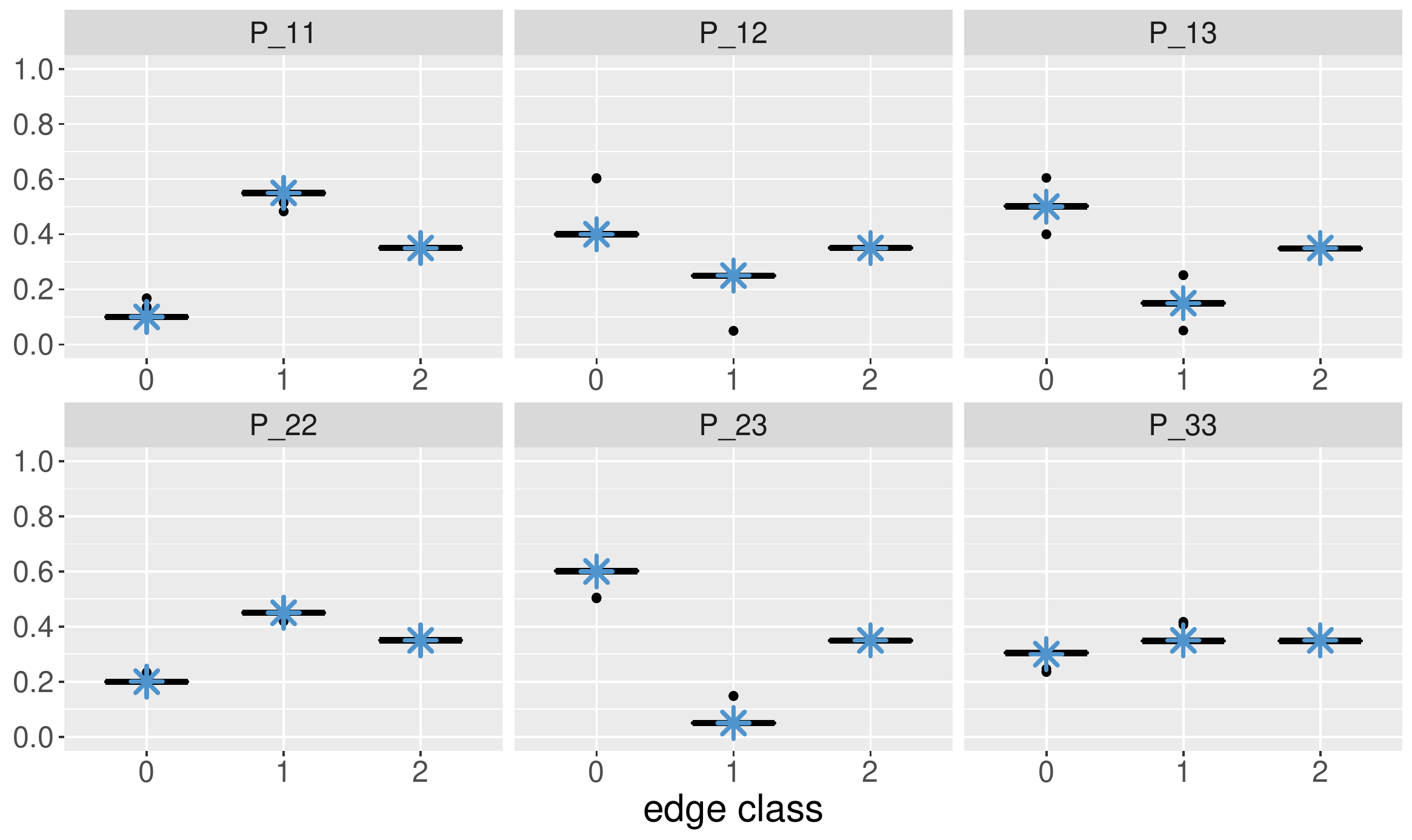}		
	}\hfill
	\caption{Boxplots of estimation results of the distributions $\bp_{ql}$ in Scenario 2 for 100 simulated networks. True underlying parameters are visualized as blue asterisks.}
	\label{fig:Sc2_Gamma}
\end{figure}

\begin{figure}
	\subfigure[$n=150$]{\includegraphics[width=0.5\textwidth]{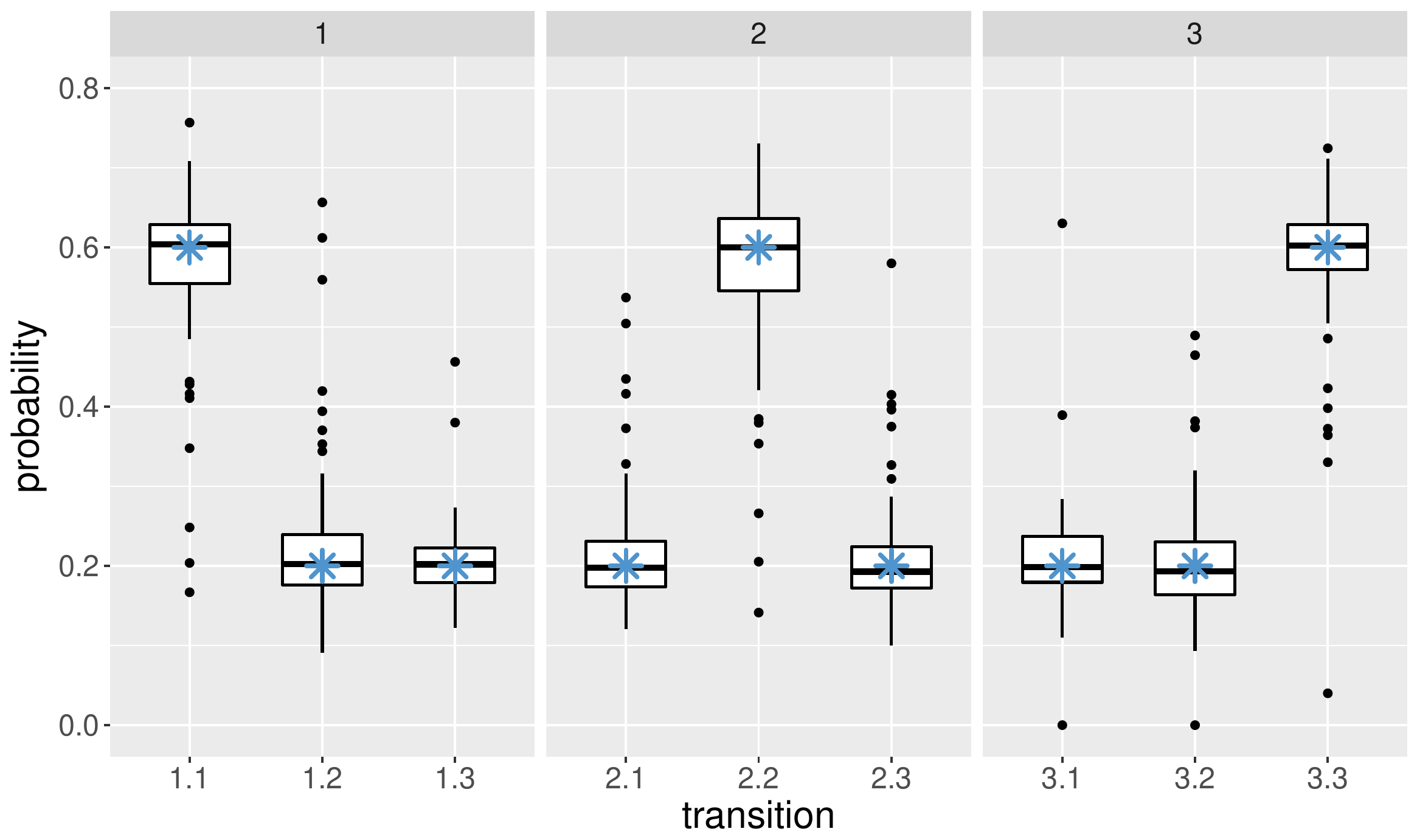}
	}\hfill
	\subfigure[$n=1000$]{\includegraphics[width=0.5\textwidth]{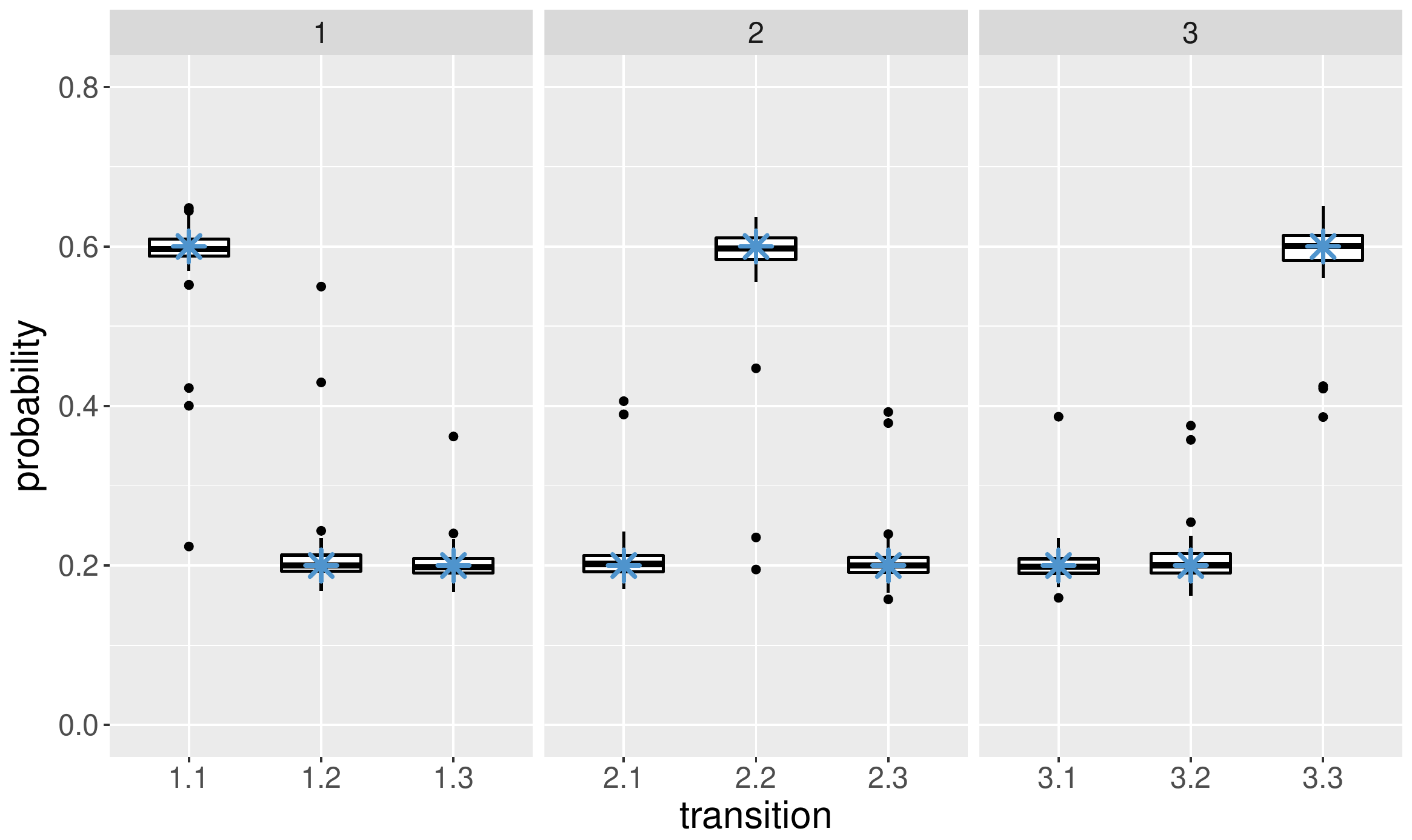}		
	}\hfill
	\caption{Boxplots of estimation results of the transition matrix $\rho$ in Scenario 2 for 100 simulated networks. True underlying parameters are visualized as blue asterisks.}
	\label{fig:Sc2_Trans}
\end{figure}

Figures \ref{fig:Sc2_pi}--\ref{fig:Sc2_Trans} visualize the results of the estimation procedure obtained after determining the global permutation. The true parameter values are presented as blue asterisks. 

The results are qualitatively similar to those in scenario 1, but with higher variability in particular in the estimation of the transition matrix $\boldsymbol{\rho}$ for the small sized network. The runtime for 100 estimations was $\sim$20 minutes for $n=125$ and $\sim$35 hours for $n=1000$.


\section{Concluding remarks}\label{sec:conclude}

The identification results in the paper justify and extend the range of applicability of the dynamic SBM from \citet{Matias2016}.
Finitely weighted edges are particularly relevant for applications, and we cover the case of a small set of edge states and a comparatively large number of latent node groups. In case of continuous edge states, it would be interesting to include a smoothing step into the variational EM algorithm, see e.g.~\citet{gassiat} in case of the EM algorithm for nonparametric hidden Markov models. 

The affiliation case, where the probabilities of the edge states only depend on whether both nodes constituting the edge are in the same state or not, no matter which particular state, still deserves further discussion. The transition matrices in the counterexample in \citet{Matias2016} are not ergodic, so that it might suffice to achieve identification with only small additional assumptions.  

It would also be of interest to extend the model to include covariates. In most applications these are node-specific, for example ethnicity in a high-school friendship network or keywords in a citation network \citep{huang2018pairwise}. For a static, binary stochastic block model with logistic link function, and where the effect of the covariates does not depend on the latent states, \citet{roy2016likelihood} discuss likelihood inference based on a Monte Carlo EM type algorithm. Our results could be directly applied if the support of the covariate vector includes a sufficiently large set on the diagonal, where they are restricted to be equal for all nodes. To take advantage of additional information of covariates for identification as e.g.~in \citet{hunter2012semiparametric} for mixture of regression models, new arguments would be required, however.

\section*{Acknowledgements}

The authors are grateful to two anonymous reviewers and to the associate editor Prof. Klaus-Robert Müller for helpful suggestions, and to Catherine Matias for comments on an earlier version of the manuscript. 

\section{Proofs}\label{sec:proofs}

\subsection{Proof of Theorem \ref{th:nonparident} and of Remark \ref{rem:inhomogeneous}}
\small{

\begin{proof}[{\sl Proof of Theorem \ref{th:nonparident}}]

As argued in the discussion preceding the theorem, under the assumption $a.$, by \citet[theorem 15]{Allman2011} the marginal parameters $\boldsymbol{\pi}$, $p_{ql}^t$ and  $P_{ql}^t$ are identified from the distribution of $\boldsymbol{X}^t$ up to permutations $\sigma_t$ of the node states $\{ 1, \ldots, Q\}$ which, however, depend on $t$. Assumption $b.$~then allows to align the node states globally, i.e.~to pass to a global permutation $\sigma$.

Then it remains to identify the transition matrix $\boldsymbol{\rho}$. To this end, fix one edge $\{i,j\}$ and a time point $t_0 \leq T-1$, and consider the joint distribution of two consecutive edge variables
	\begin{align*}
	\nonumber
	\mathbb{P}(X_{ij}^{t_0}, X_{ij}^{t_0+1})&=\sum_{q_1,l_1}\sum_{q_2,l_2}\mathbb{P}(X_{ij}^{t_0}, X_{ij}^{t_0+1} \cond Z_i^{t_0}=q_1,Z_j^{t_0}=l_1,Z_i^{t_0+1}=q_2, Z_j^{t_0+1}=l_2)\\
	\nonumber
	&\qquad\qquad\qquad\cdot \mathbb{P}(Z_i^{t_0}=q_1,Z_j^{t_0}=l_1,Z_i^{t_0+1}=q_2, Z_j^{t_0+1}=l_2)\\ \nonumber
	&=\sum_{q_1,l_1}\sum_{q_2,l_2} \mathbb{P}(X_{ij}^{t_0}\cond Z_i^{t_0}=q_1, Z_j^{t_0}=l_1)\mathbb{P}(X_{ij}^{t_0+1}\cond Z_i^{t_0+1}=q_2, Z_j^{t_0+1}=l_2)\\ \nonumber
	&\qquad\qquad\qquad\cdot \mathbb{P}(Z_i^{t_0}=q_1,Z_i^{t_0+1}=q_2)\mathbb{P}(Z_j^{t_0}=l_1,Z_j^{t_0+1}=l_2)\\
	&=\sum_{q_1,l_1}\sum_{q_2,l_2} \pi_{q_1}\pi_{l_1}\rho_{q_1q_2}\rho_{l_1l_2}~\big(\mu_{q_1l_1}^{t_0}\times \mu_{q_2l_2}^{t_0+1}\big).
	\label{dynamischbivariat}
	\end{align*}
	Using that by assumption, $\mu_{ql}^t=\mu_{lq}^t$ for every $t$, by collecting terms we may write the last expression as 
\begin{equation}
	\label{FFF}
	\mathbb{P}(X_{ij}^{t_0}, X_{ij}^{t_0+1})=\sum_{q_1\leq l_1}\sum_{q_2\leq l_2} \varphi_{q_1,q_2,l_1,l_2}~~\big(\mu_{q_1l_1}^{t_0}\times \mu_{q_2l_2}^{t_0+1}\big),
	\end{equation}
where
	$$ \varphi_{q_1,l_1,q_2,l_2}=
	\begin{cases}
	\pi_{q_1}^2\rho_{q_1q_2}^2, &\quad \text{if } q_1=l_1,q_2=l_2,\\
	2\pi_{q_1}^2\rho_{q_1q_2}\rho_{q_1l_2}, &\quad \text{if } q_1=l_1,q_2<l_2, \\
	2\pi_{q_1}\pi_{l_1}\rho_{q_1q_2}\rho_{l_1q_2}, &\quad \text{if } q_1<l_1,q_2=l_2, \\
	2\pi_{q_1}\pi_{l_1}(\rho_{q_1q_2}\rho_{l_1l_2} + \rho_{q_1l_2}\rho_{l_1q_2}), &\quad \text{if } q_1<l_1,q_2<l_2. 
	\end{cases}
	$$

	We denote by $\varPhi$ the ${Q+1\choose 2} \times {Q+1 \choose 2}$ square matrix with entries 
	$$\varPhi\big(\small{(q_1,l_1);(q_2,l_2)}\big)=\varphi_{q_1,l_1,q_2,l_2}$$
	with rows indexed by $(q_1,l_1)$ and columns indexed by $(q_2,l_2)$.
	
	Let $G_{ql}^t$ denote the cumulative distribution function of $\mu_{ql}^t$.
	Since measures $\mu_{ql}^t$, $1 \leq q \leq l \leq Q$ in (\ref{eq:conddistr}) are linearly independent by assumption $a.$, from \citet[lemma 17]{Allman2009} it follows that for each $t$ there exist ${ Q+1 \choose 2}$ points $\{u_{ql}^t\}$ such that the ${Q+1 \choose 2} \times { Q+1 \choose 2}$ matrices  
	$$M_t\big(\small{(q_1,l_1);(q_2,l_2)}\big)=G_{q_1l_1}^t(u_{q_2l_2}^t), \qquad 1 \leq q_i \leq l_i \leq Q \text { and } i=1,2,$$ 
	have full rank. 

	Further, we let $F$ denote the distribution function of (\ref{FFF}) and let $N_{t_0}$ be the ${Q+1 \choose 2} \times { Q+1 \choose 2}$ sized matrix with entries $$N_{t_0}\big(\small{(q_1,l_1);(q_2,l_2)}\big)=F(u_{q_1l_1}^{t_0},u_{q_2l_2}^{t_0+1}).$$ Then from (\ref{FFF}) it follows that 
	$N_{t_0}=M_{t_0}' \varPhi M_{t_0+1}$ and hence 
	$$\varPhi=(M_{t_0}')^{-1}N_{t_0}M_{t_0+1}^{-1},$$
	where the right-hand side identifies the left-hand side. Since the entries of $\boldsymbol{\pi}$ are also identified, we obtain identification of $\boldsymbol{\rho}$ from
	 $$ \rho_{ql} = \sqrt{\varphi_{q,q,l,l}}/\pi_q.$$
\end{proof}

\begin{proof}[{\sl Proof of statement in Remark \ref{rem:inhomogeneous}}]
	The marginal  distribution of  $(\boldsymbol{Z}^1,\boldsymbol{X}^1)$ is in $\mathcal{M}(n,Q, \mathcal{X}, \cP)$ with parameters $\boldsymbol{\pi} = \boldsymbol{\pi}^1$, $p_{ql}^1$ and  $P_{ql}^1$, while the marginal  distributions of  $(\boldsymbol{Z}^t,\boldsymbol{X}^t)$ have parameters $p_{ql}^t$ and  $P_{ql}^t$ and $\boldsymbol{\pi}^t = \boldsymbol{\pi} \cdot \boldsymbol{\rho^2} \cdot \ldots \cdot \boldsymbol{\rho}^t$, $t=2, \ldots, T$. By irreducibility of the $\boldsymbol{\rho}^s$, $s=2, \ldots, t$, the entries of the $\boldsymbol{\pi}^t$ are strictly positive.
	By assumption $a.$ and since $n \geq 9$, from \citet[theorem 15]{Allman2011} the marginal parameters $\boldsymbol{\pi}$, $p_{ql}^t$, $P_{ql}^t$ and also $\boldsymbol{\pi}^t$ are identified from the distribution of $\boldsymbol{X}^t$ up to permutations $\sigma_t$ of the node states $\{ 1, \ldots, Q\}$ which, however, depend on $t$. Assumption $b.$~again allows to align the node states globally, i.e.~to pass to a global permutation $\sigma$.
	 In the argument of the proof of Theorem \ref{th:nonparident}, the transition matrix $\boldsymbol{\rho}^t$, $t=2, \ldots, T$ is identified based on the joint distribution of $(X_e^{t-1}, X_e^{t})$ for some edge $e$, and knowledge of $\boldsymbol{\pi}^{t-1}$ and of $\mu_{ql}^{t-1}$, $\mu_{ql}^{t}$, $1 \leq q \leq l \leq Q$.   
	Hence the argument also works in the inhomogeneous case.
\end{proof}

}

\subsection{Proof of Theorem \ref{th:smallstatic}}
{\small
First, note that all vectors $\boldsymbol{p}_{ql}$ must be pairwise different, for, if two vectors were equal, $\cC_{m,Q,\kappa}$ would not have full row rank. 
Now, as $\cC_{m,Q,\kappa}$ has full row rank, \citet[lemma 17]{Allman2011} provides three pairwise disjoint subsets $G_i \subset E$ of edges, $i=1,2,3$  each consisting of the union of $m$ complete graphs over $m$ different nodes and hence of cardinality $|G_i| = m {m \choose 2}$, such that the $Q^n \times \kappa^{|G_i|}$ matrices $\bar \cC^{(i)}$ which contain the distribution of $(X_e)_{e \in G_i}$ conditional on node assignments have full row rank.  
Thus, by \citet[theorem 4a]{Kruskal}, see also \citet[theorem 3]{Rhodes2010},  in the version of \citet[theorem 16]{Allman2011} we can identify each of the matrices $\bar \cC^{(i)}, i=1,2,3$ as well as the vector $\boldsymbol{\lambda} \in [0,1]^{Q^n}$ containing the probabilities for each node state assignment up to a permutation of the rows. The entries of $\boldsymbol{\lambda}$ are of the form $\prod_{q=1}^Q \pi_q^{m_q}$ for $m_q$ denoting the number of nodes which are in state $q$.

To complete the proof, we must identify the vectors $\fpi$ and $\boldsymbol{p}_{ql}, 1 \leq q \leq l \leq Q.$
We consider one of the complete subgraphs on $m$ nodes in $G_1$, which we may assume to be the nodes $\{1, \ldots, m\}$. Its conditional distribution matrix given node states is $\cC_{m,Q,\kappa}$, which is identified up to row permutation. 

We are able to recover the whole set of parameters $\{\bp_{ql}\}$ (without assignment) by doing the following marginalizations.
We fix an edge and for each row  $\bx$ of (unknown) node assignments, we perform marginalizations by summing all columns in which the specific edge is in state $k=1,\ldots,\kappa$. Now, we put these values into a vector of size $\kappa$. The resulting vector is equal to a parameter vector $\bp_{ql}$ but we do not know $q$ and $l$ so far. To assign these, we consider for every row of $\cC_{m,Q,\kappa}$ the set $A_{\bx}$, which contains all the vectors we obtained by marginalizations in this row. This set is of the form 
$$A_{\bx}:=\{\bp_{ql}\,|\, \text{two nodes are in state $q$ and $l$ in assignment $\bx$}\}.$$

If all node states in assignment $\bx$ are equal, then for each edge, we recover the same vector $\bp_{qq}, q =1,\ldots, Q$ and this will be the only element contained in $A_{\bx}$.
In contrast, if at least one node is in a distinct state, since $m \geq 3$ there are assignments which result in two distinct conditional edge distributions $\bp_{ql}$ and $\bp_{\tilde q \tilde l}$, where edge states $q,l,\tilde q, \tilde l$ are all contained in $\bx$ (some may be equal). Since $\bp_{ql} \not= \bp_{\tilde q \tilde l}$, the set $A_{\bx}$ will contain more than one element.
Hence, there are exactly $Q$ rows for which $A_{\bx}$ is a one element set. Choosing an arbitrary labeling, we have thus recovered $\bp_{qq}, 1 \leq q \leq Q$.
Futher, the entries of $\boldsymbol{\lambda}$ belonging to these rows are of the form $\pi_q^m$. Thus, extracting the $m$-th roots and combining the values to a vector gives us $\fpi$. 

To finally assign the vectors $\bp_{ql}, q\neq l$, we observe that the rows for which $A_{\bx}$ has precisely two elements are exactly those rows for which all but one of the nodes are in the same state $q$ and one is in another state $l\neq q$, and then $A_{\bx}=\{\bp_{qq},\bp_{ql}\}$. 
For, if there is at least one additional node in a distinct state $\tilde l\neq q$ (which may also be equal to $l$), since $m \geq 3$, there are assignments resulting in $\bp_{l\tilde l}$ which is different from $\bp_{qq}$ and $\bp_{ql}$. So, $A_{\bx}$ contains more than two elements in all other cases.

The missing parameters $\bp_{ql}$ are contained in sets either together with $\bp_{qq}$ or together with $\bp_{ll}$. As these already are identified, we can match the sets and get the assignments for $\bp_{ql}, q \neq l$.
The result follows.
}
\subsection{Proof of Theorem \ref{th:dynfewedgestates}}
{\small

We shall use the following additional notation.  If $M_j$ are matrices of dimension $r \times k_j$, $j=1,2,3$, then the matrix triple product of the $M_j$ is the three-fold tensor of dimension $k_1 \times k_2 \times k_3$ defined by
\[ [M_1,M_2,M_3](i_1,i_2,i_3) = \sum_{l=1}^r \, (M_1)_{l,i_1} \, (M_2)_{l,i_2} \, (M_3)_{l,i_3}, \quad 1 \leq i_j \leq k_j,\ j=1,2,3. \]

If, in addition, $v$ is a row vector of dimension $r$, and $\text{diag}(v)$ the $r \times r$ diagonal matrix with entries $v$ on its diagonal, then we set
\[ [v;M_1,M_2,M_3] = [\text{diag}(v) \cdot M_1,M_2,M_3].\]

Turning to the proof of Theorem \ref{th:dynfewedgestates}, as argued in the discussion preceding the theorem, from Theorem \ref{th:smallstatic} the marginal parameters $\boldsymbol{\pi}$ and $p_{ql}^t$ are identified from the distribution of $\boldsymbol{X}^t$ up to permutations $\sigma_t$ of the node states $\{ 1, \ldots, Q\}$. By assumption, these may then be aligned globally.

In order to identify the transition matrix $\f{\rho}$, we shall apply Kruskal's theorem as in \citet[theorem 16]{Allman2011}, and note that the submodel on $m$ nodes, $(\boldsymbol{X}^t_m)_{t=1\ldots, T}$ is an HMM with 
\[ \text{transition matrix } \f{\varrho}:=\f{\rho}^{\otimes m}, \qquad \text{stationary distribution } \f{\varpi}:=\f{\pi}^{\otimes m}.\]
	We fix some time point $2 \leq t_0 \leq T-1$. As in previous arguments on nonparametric HMM identifiability \citep{Allman2009, gassiat, Holzmann2016}, conditionally on $(Z^{t_0}_1, \ldots, Z^{t_0}_m)'$ the random vectors $\f{X}^{t_0-1}_m,\f{X}^{t_0}_m,\f{X}^{t_0+1}_m$ are independent. The conditional distribution of $\f{X}_m^{t_0}$ is given by $\mathcal{C}^{t_0}_{m,Q,\kappa} = \mathcal{C}_m^{t_0}$ as defined in (\ref{eq:thecondmatrix1}), where we drop the subscripts $Q$ and $\kappa$ in the following.  
	Then, the distribution of $\f{X}_m^{t_0+1}$ conditional on $(Z^{t_0}_1, \ldots, Z^{t_0}_m)'$ is given by $\f{\varrho}\,\mathcal{C}_m^{t_0+1}$, and that of $\bX_m^{t_0-1}$ by $\widetilde{\f{\varrho}}\mathcal{C}_m^{t_0-1}$, where $\widetilde{\f{\varrho}}$ is the time reversal of $\f{\varrho}$,
\[ 	\widetilde{\f{\varrho}}:=\diag(\f{\varpi}^{-1}) \f{\varrho}' \diag(\f{\varpi}) = \widetilde{\f{\rho}}^{\otimes m}, \qquad \widetilde{\f{\rho}}:=\diag(\f{\pi}^{-1}) \f{\rho}' \diag(\f{\pi}).\]
The joint distribution can therefore be written as \[ \mathbb{P}(\f{X}_m^{t_0-1},\f{X}_m^{t_0},\f{X}_m^{t_0+1})=[\f{\varpi};\f{\varrho}'\mathcal{C}_m^{t_0-1},\mathcal{C}_m^{t_0},\f{\varrho}\,\mathcal{C}_m^{t_0+1}].\]
	By assumption $\mathcal{C}_m^{t}$ has full row rank for every $t$. Further, $\f{\rho}$ has full rank by assumption, and hence so does its $m$-fold tensor product $\f{\varrho}=\f{\rho}^{\otimes m}$. Moreover, since $\f{\pi}$ has no zero entries since $\f{\rho}$ is ergodic, also the time reversals $\widetilde{\f{\rho}}$ as well as $\widetilde{\f{\varrho}}$ have full rank.
	 Therefore, it follows that $\widetilde{\f{\varrho}} \mathcal{C}_m^{t_0-1}$ and $\f{\varrho}\,\mathcal{C}_m^{t_0+1}$ also have full row rank. 
	Since $\f{\varpi}$ does not have zero entries, from Kruskal's theorem \citet{Kruskal} in the version of \citet[theorem 16]{Allman2011} we in particular identify the matrices $\f{\varrho}\,\mathcal{C}_m^{t_0+1}$ and $\mathcal{C}_m^{t_0}$ up to permutation, and the matrix $\mathcal{C}_m^{t_0}$ allows to align permuations.  Since $\mathcal{C}_m^{t_0+1}$ also has previously been identified, we identify $\f{\varrho}$ by
	\[ \f{\varrho} = \f{\varrho}\,\mathcal{C}_m^{t_0+1}\, (\mathcal{C}_m^{t_0+1})'\, \big(\mathcal{C}_m^{t_0+1}\, (\mathcal{C}_m^{t_0+1})'\big)^{-1}.\]
	The parameters of $\f{\rho}$ are extracted as the $m$-th root of particular matrix entries of $\f{\varrho}$, namely for $\bz_q=(q,\ldots,q)$ and $\bz_l=(l,\ldots,l)$ it holds that
	$$\rho_{ql}=\big(\f{\varrho}(\bz_q,\bz_l)\big)^{1/m}.$$

}

\begin{proof}[{\sl Proof of statement in Remark \ref{rem:inhomogeneous2}}]
	Arguing as before the marginal parameters $\boldsymbol{\pi}^t$ and $\bp_{ql}^t$, $t=1, \ldots, T$, are identified, and it remains  
	to identify the transition matrices $\boldsymbol{\rho}^t$, $t=2, \ldots, T$. As in the above argument with consider the (inhomogenous) HMM $(\boldsymbol{X}^t_m)_{t=1\ldots, T}$   with 
	\[ \text{transition matrices } \f{\varrho}^t:=(\f{\rho}^t)^{\otimes m}, \qquad \text{marginal state distributions } \f{\varpi}^t:=(\f{\pi}^t)^{\otimes m}.\]
	 Given a time point $2 \leq t \leq T-1$, conditionally on $(Z^{t_0}_1, \ldots, Z^{t_0}_m)'$ the random vectors $\f{X}^{t-1}_m,\f{X}^{t}_m,\f{X}^{t+1}_m$ are independent with 
	 $\mathcal{C}^{t}_{m,Q,\kappa} = \mathcal{C}_m^{t}$ as conditional distribution of $\f{X}_m^{t}$, $\f{\varrho}^{t+1}\,\mathcal{C}_m^{t+1}$ as conditional distribution of $\f{X}_m^{t+1}$ and $\widetilde{\f{\varrho}}^{t}\,\mathcal{C}_m^{t-1}$ as conditional distribution of $\f{X}_m^{t-1}$
	 where 
	\[ \widetilde{\f{\varrho}}^t = \diag(\f{\varpi}^{t})^{-1}\, (\f{\varrho}^t)'\, \diag(\f{\varpi}^{t-1}).\]   
	By our full-rank assumption, as above we identify $\f{\varrho}^3, \ldots, \f{\varrho}^T$ and hence $\boldsymbol{\rho}^3, \ldots, \boldsymbol{\rho}^T$. Moreover, we also identify $\widetilde{\f{\varrho}}^2$, and since $\diag(\f{\varpi}^{s})$, $s=1,2$ are also identified, we obtain $\f{\varrho}^2$ and hence $\boldsymbol{\rho}^2$.  
\end{proof}

\bibliographystyle{chicago}

\bibliography{Literatur}

\end{document}